\documentclass[twoside,11pt]{amsart}
\usepackage{graphicx, lscape}
 \usepackage{amsaddr}
\usepackage{amsmath, amsfonts, amssymb, ifthen}
\usepackage{caption}
\usepackage{subcaption}
\usepackage[justification=centering]{caption}
\usepackage[colorlinks=true]{hyperref}% this thing makes links like thereoms and citations incolors

\newcommand{\CC}{\mathbb{C}}
\newcommand{\DD}{\mathbb{D}}
\newcommand{\RR}{\mathbb{R}}
\newcommand{\EE}{\mathcal{E}}
\newcommand{\mandel}{\mathcal{M}}
\newcommand{\crpt}{\xi}

\newcommand{\Rnca}{R_{n,a,c}}
\newcommand{\KRnca}{K(R_{n,a,c})}
\newcommand{\Rnta}{R_{n,a,ta}}
\newcommand{\Mn}{{M}_n}
\newcommand{\Wcj}{W_{c,j}}
\newcommand{\Vcj}{V_{c,j}}
\newtheorem{theorem}{Theorem}[section]

\newtheorem{main}{Main Theorem}

\theoremstyle{definition}

\newtheorem*{notation}{Notation}

\newcommand{\ep}{\varepsilon}

\DeclareMathOperator\Arg{Arg}

\newtheorem{prop}[theorem]{Proposition}
\newtheorem{cor}[theorem]{Corollary}
\newtheorem{lem}[theorem]{Lemma}
\newtheorem{defn}[theorem]{Definition}

\begin{document}

\title{Baby Mandelbrot sets and Spines in some one-dimensional subspaces of the parameter space for generalized McMullen Maps}

\author[S.~Boyd]{Suzanne Boyd}
\address{Department of Mathematical Sciences\\
University of Wisconsin Milwaukee\\
PO Box 413\\
Milwaukee, WI 53201, 
USA}
\email{sboyd@uwm.edu, ORCID: 0000-0002-9480-4848}

\author[M.~Hoeppener]{Matthew Hoeppner}
\address{Mathematics and Statistics Department\\
Hope College\\
27 Graves Place\\
Holland, MI 49423, 
USA}
%\email{hoeppner@hope.edu}

\date{\today}

\begin{abstract}

For the family of complex rational functions of the form 

\noindent \mbox{$\Rnca(z) = z^n + \dfrac{a}{z^n}+c$,} 
 known as ``Generalized McMullen maps'', for $a\neq 0$ and $n \geq 3$ fixed, we study the boundedness locus in some one-dimensional slices of the $(a,c)$-parameter space, by fixing a parameter or imposing a relation. 
 
 First, if we fix  $c$ with $|c|\geq 6$ while allowing $a$ to vary, assuming a modest lower bound on $n$ in terms of $|c|$, we establish the location in the $a$-plane of $n$ ``baby" Mandelbrot sets, that is, homeomorphic copies of the original Mandelbrot set. We use polynomial-like maps, introduced by Douady and Hubbard (\cite{DouadyHubbard}) and applied for the subfamily \mbox{$R_{n,a,0}$} by Devaney (\cite{DevaneyHalos}).  

Second, for slices in which $c=ta$, 
we again observe what look like baby Mandelbrot sets within these slices, and begin the study of this subfamily by establishing a neighborhood containing the boundedness locus.
\end{abstract}

\maketitle

\markboth{\textsc{S. Boyd and M. Hoeppner}}
  {\textit{}}


\footnotetext[2 ]{All authors contributed equally to this paper. We would like to thank Brian Boyd for the computer program detool which we used to generate all of the Mandelbrot and Julia images in this paper.}

%%%%%%%%%%%%%%%%%%%%%%%%%%%%%%%%%%%%%%%%%%%%%%%%
\section{Introduction} 
\label{sec:introduction}
%%%%%%%%%%%%%%%%%%%%%%%%%%%%%%%%%%%%%%%%%%%%%%%%

As a simple starting example, consider the family of quadratic polynomials
$
P_c(z) = z^2 + c,~c \in \mathbb{C}.
$
We define the \textit{Fatou set} of $P_c$ in the typical way, as the set of values in the domain where the iterates of $P_c$ is a normal family in the sense of Montel. The \textit{Julia set}, $J$, is also defined the usual way as the complement to the Fatou set. The \textit{filled Julia set}, $K$, is the union of the Julia set and the bounded Fatou components.

\begin{figure}
\centering
\includegraphics[width=.45\textwidth,keepaspectratio]{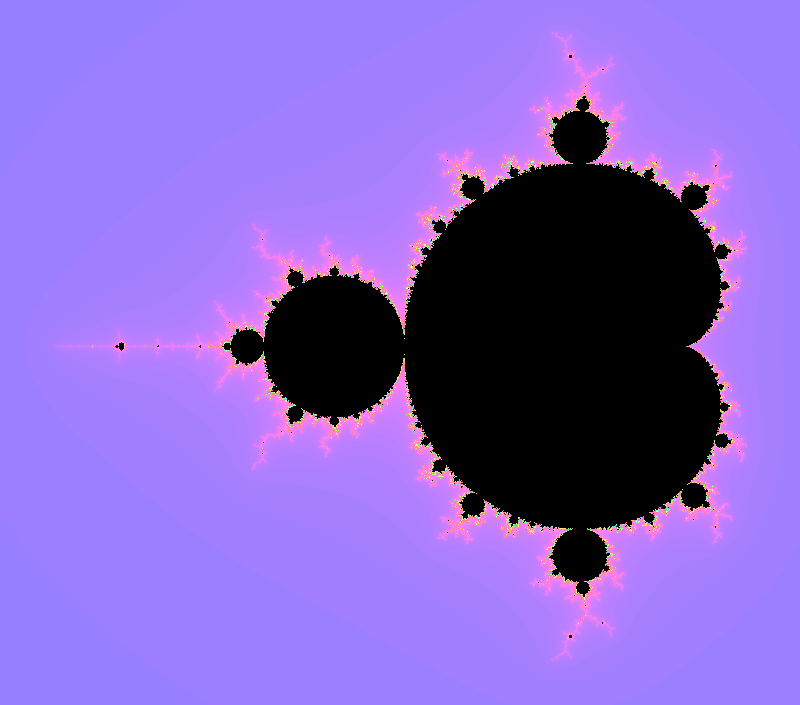}
\includegraphics[width=.4\textwidth,keepaspectratio]{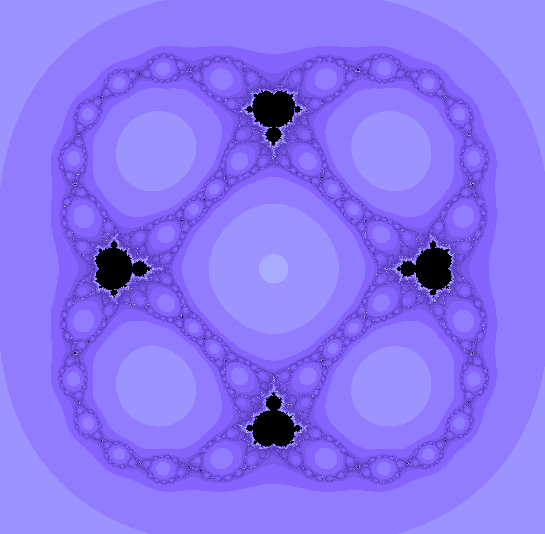}
\caption{Parameter Planes of $P_c$ (left) and $R_{5,a,0}$ (right).}
\label{fig:Mandelbrotone}
\end{figure}

The \textit{Mandelbrot Set}, $\mandel$, is the set of $c$-values such that the critical orbit of $P_c$ is bounded, here that is the orbit of $0$.  Figure \ref{fig:Mandelbrotone} (left) is the Mandelbrot set drawn in the $c$-parameter plane of $P_c$.  For other functions, the set of parameter values where at least one critical orbit is bounded is called the \textit{boundedness locus}.

The study of the Mandelbrot set became more accessible as computers advanced.  Douady and Hubbard (\cite{DouadyHubbard}) showed this set can result from other iterative processes as well, showing multiple homeomorphic copies of the Mandelbrot set occur when Newton's Method is applied to a cubic polynomial family with a single parameter. They defined what it means for a map to behave like $P_c$, calling such a map {\em polynomial-like of degree two} (see Section \ref{sec:prelim}). Mc{M}ullen (\cite{mcmullen}) shows that every non-empty bifurcation locus of any analytic family will contain quasiconformal copies of the Mandelbrot set of $P_c$ (or of $z\mapsto z^n + c$, based on the multiplicity of critical points), but we are interested in using Douady and Hubbard's approach to prove that Mandelbrot set copies exist in some specific locations in some parameter planes, for the  family:
$$
\Rnca(z) = z^n + \dfrac{a}{z^n}+c~,~n \in \mathbb{N}, n \geq 3, ~a \in \mathbb{C}^*,~c \in \mathbb{C},
$$
where $\mathbb{C}^* = \mathbb{C}\backslash \lbrace 0\rbrace$. We consider $n$ a fixed integer with $n \geq 3$, thus have parameter space $\CC^2$.
We let $\Mn(\Rnca)$ represent the boundedness locus of $\Rnca$ for a fixed $n$. We call homeomorphic copies of the Mandelbrot set, $\mandel$, ``baby'' Mandelbrot sets or ``baby $\mandel$'s''.

This family, including the one-parameter subfamily with $c=0$, has been studied previously by Devaney and colleagues as well as the first author and colleagues.
In \cite{DevaneyRussell,DevaneyKozma}, Devaney and coauthors study the family in the case of $c$ at the center of a hyperbolic component of the Mandelbrot set for $P_c$ (that is, the critical point is a fixed point). 
 For $n\geq 2$, Devaney and colleagues study the subfamily with $c=0$, ``McMullen maps", in papers such as \cite{DevaneyHalos,DevaneySurvey2013}. 
 They establish the location of $n-1$ baby $\mandel$'s in the slice of $\Mn(\Rnca)$ in the $a$-parameter plane when $c=0$ (see Figure \ref{fig:Mandelbrotone} (right) for an example). 
 In \cite{devgar}, Devaney and Garijo  study Julia sets as the parameter $a$ tends to $0$, for a different generalization of McMullen Maps: $z \mapsto z^n + \dfrac{a}{z^d}$ with $n \neq d$.
In \cite{DevaneyHalos,jangso} the authors find $n$ baby $\mandel$'s in the $n\neq d$ case.

In \cite{BoydSchulz}, the first author and Schulz study the geometric limit as $n\to \infty$ of Julia sets and of $\Mn(\Rnca)$, for $\Rnca$ for any complex $c$ and any complex $a\neq 0$.  
In \cite{BoydMitchell}, the first author and Mitchell establish the location of baby $\mandel$'s in both the $a$ and $c$-parameter planes of $\Rnca$, in the case $c\neq 0$ (and $a\neq 0$). First, in the $a$-plane they study the case $c\in[-1,1]$, establishing a baby $\mandel$ in the case of $n\geq 3$ and $c\in [-1,0],$ and in the case of odd $n\geq 3$ and $c\in [0,1].$  Second, in the $c$-plane, for $n\geq 5$ and $1 \leq a \leq 4$, they locate $n$ baby $\mandel$'s, and $2n$ if $n$ is odd; additionally, they locate a baby $\mandel$'s in the $c$-plane for $n \geq 11$ and $\frac{1}{10} \leq a \leq 1$, and show if $n$ is odd there are at least two.  Note in both of these results, the fixed parameter is always real.

Not all of the behavior of this family is polynomial like. For example, there are Julia sets which are Cantor sets of simple closed curves (\cite{mcmullen_example, devlook}). 
 Xiao, Qiu, and Yongchen (\cite{xiaoqiu}) establish a topological description of the Julia sets (and Fatou components) of $\Rnca$ according to the dynamical behavior of the orbits of its free critical points. This work includes a result that if there is a critical component of the filled Julia set which is periodic while the other critical orbit escapes, then the Julia set consists of infinitely many homeomorphic copies of a quadratic Julia set, and uncountably many points. However in order to find baby $\mandel$'s in our parameter planes of interest, we first locate baby Julia sets,  using different techniques (based on specific parameter ranges rather than the type of dynamical behavior). 

%%%%%%%%%%%%%%%%%%%%%%%%%%%%%%%%%%%%%%%%%%%%%%%%%%%%%%

This family has $2n$ critical points. While this seems like a lot of critical orbits to keep track of, in fact every one of these critical points maps to one of just two critical values, $v_{\pm} = c\pm2\sqrt{a}$, where we choose the branch cut of the square root function to be the negative real axis. 
We also have two free complex parameters $a$ and $c$. So we are essentially studying the dynamics in the case of two degrees of freedom; this is inherently a two complex dimensional parameter space. 
There are various approaches to the study of a two complex dimensional parameter space.

In this article, we study $M_n(\Rnca)$ in 
one-dimensional subfamilies or ``slices'' of parameter space. One type of slice is to fix one parameter $a$ or $c$, and let the other vary (horizontal and vertical lines in the parameter space). A second type of slice could be created by allowing both parameters $a$ and $c$ to vary in a one-dimensionally constrained way; in our example, we consider the lines $c=ta$ for any non-zero complex slope $t$ (lines through the origin in the parameter space). 

In any of these one-complex-dimensional subfamilies, the goal is the same: establish locations of baby $\mandel$'s for specific parameters. 
Our strategy is to go about achieving this goal in a slice of interest by completing some or all of the following steps, not necessarily in this order:
\begin{enumerate}
    \item Identify a broad region (such as a large disk or polar rectangle) in which the $M_n(\Rnca)$ must lie in the given slice.
\item Determine a curve, or ``spine'', around which $M_n(\Rnca)$ is centered in the slice. Ideally, find a neighborhood of the spine which which contains $M_n(\Rnca)$ in the slice. 
\item Identify the collection of parameter values in the given slice for which the critical points are fixed points. This helps to find the centers of potential baby $\mandel$'s, which serve as anchors for the creation of sets to contain each individual baby $\mandel$.
\item Identify and prove a location of one or more baby $\mandel$s.
\end{enumerate}

For parameter slices with $a$ or $c$ fixed, \cite{BoydSchulz} gets one to step 2.  Then, in \cite{BoydMitchell} the location of baby $\mandel$'s in the case of $c\in[-1,1]$ in the $a$-plane is determined.
Similarly, in the first result of this article, we build on this work to locate baby $\mandel$'s in the $a$-plane, but this time for $c$ larger and complex. 

\begin{main}
\label{thm:Main_Theorem_APlane}
For $|c| \geq 6$ and $n$ such that $4|c| + 8 \leq 2^{n+1}$, there exist $n$ ``primary'' homeomorphic copies of the Mandelbrot set in the boundedness locus, $M_n(\Rnca)$, in the $a$-plane of the family  $\Rnca(z) = z^n + \frac{a}{z^n} + c$.
\end{main}

\begin{figure}
\centering
\begin{subfigure}{.5\textwidth}
  \centering
  \includegraphics[width=.95\textwidth,keepaspectratio]{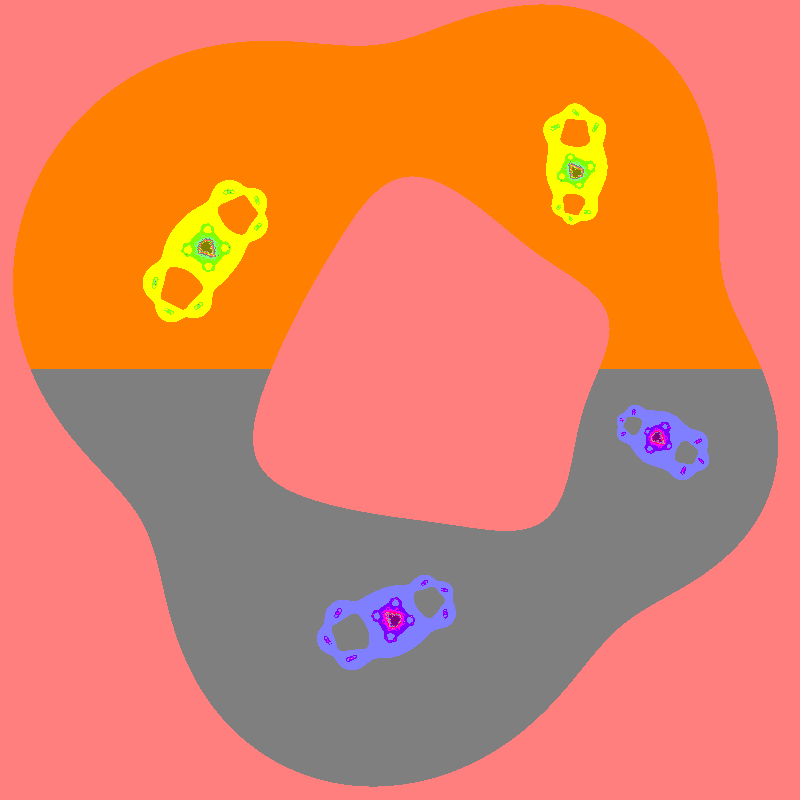} 
  \caption{The $a$-parameter plane of $R_{4,a,6i}$. }
  \label{a-plane_n4c6i}
\end{subfigure}%
\begin{subfigure}{.5\textwidth}
  \centering
 \includegraphics[width=.95\textwidth,keepaspectratio]{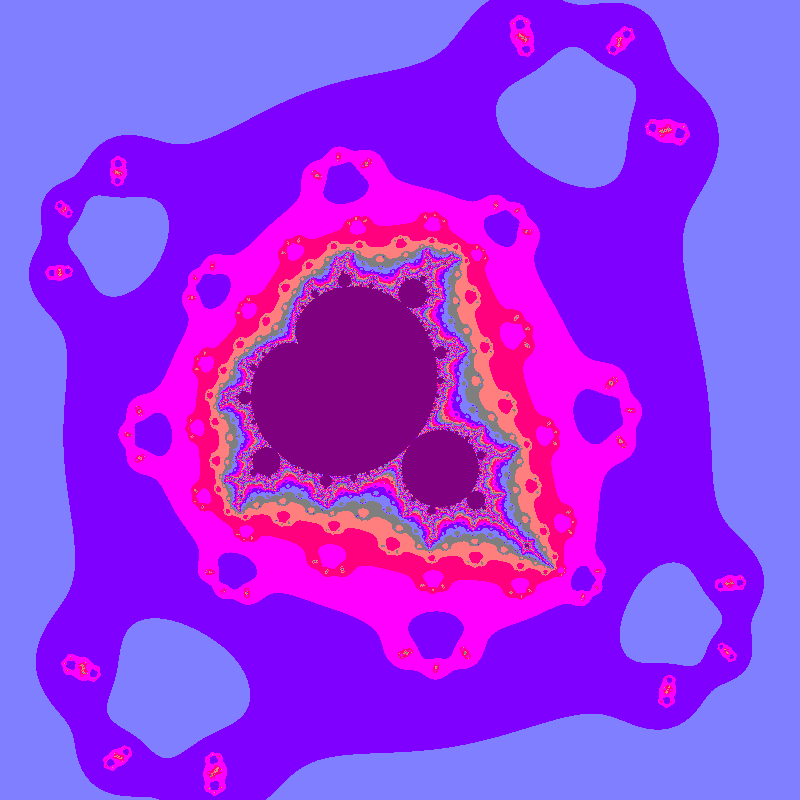} 
  \caption{A `baby' Mandelbrot set in $R_{4,a,6i}$}
  \label{a-plane_n4c6izoom}
\end{subfigure}%
\caption{The boundedness locus for $R_{4,a,6i}$. }
\label{fig:Parameter_planes_Rnca}
\end{figure}

We establish this result in Section~{\ref{sec:babyMs_Aplane}}.
See Figure~\ref{fig:Parameter_planes_Rnca}. Each pixel is assigned two RGB values, one for each of the two critical values, which are then averaged. Bounded orbits are assigned black. The color change in this figure is due to a branch cut of the square root function, since only one critical orbit is bounded in this plane. We discuss furhter how the parameter planes of $\Rnca$ are drawn in Section~\ref{sec:prelim}. 

\begin{figure}
\centering
\begin{subfigure}{.5\textwidth}
  \centering
\includegraphics[width=.95\textwidth,keepaspectratio]{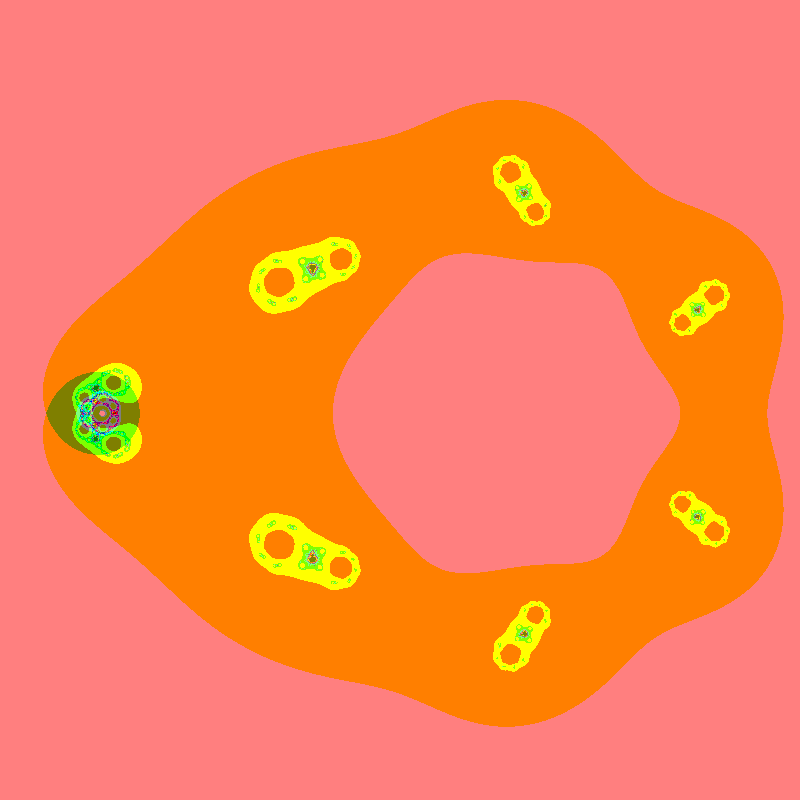}
  \caption{Parameter plane slice for $n=6$}
  \label{c=a_n=6_aplane_intro}
\end{subfigure}%
\begin{subfigure}{.5\textwidth}
  \centering
\includegraphics[width=.95\textwidth,keepaspectratio]{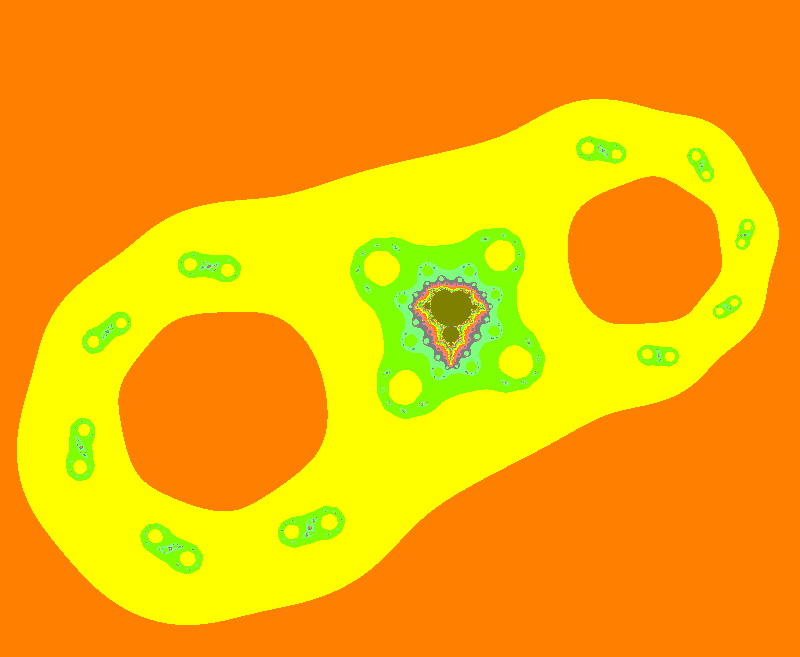}
  \caption{Zoomed in view}
  \label{c=a_n=6_mandel_intro}
\end{subfigure}%
\caption{The boundedness locus for $\Rnta(z)= z^n + \frac{a}{z^n} + ta$, for $n=6$ and $t=1$.}
\label{fig:c=a_n=6_intro}
\end{figure}

Next, we consider a different type of slice other than a fixed parameter: we set $c=ta$ for any non-zero complex slope $t$. These slices also appear to contain baby $\mandel$'s (see Figure \ref{fig:c=a_n=6_intro}), but in this case, \cite{BoydSchulz} does not provide a neighborhood for $M_n(\Rnta)$ in these types of slices. Thus we complete Steps 1 and 2 in the study of this slice, by finding a curve or \textit{spine}, $\mathcal{S}_t$, and establishing that $M_n(\Rnta)$ must lie within a neighborhood of this spine. (We take advantage of bounds on the Julia sets of these maps established in \cite{BoydSchulz}).
See Equation~\ref{eqn:spine} for the definition of $\mathcal{S}_t$.

% SECOND MAIN THEOREM
\begin{main}
\label{thm:Main_Theorem_C=ta}
Let $\epsilon > 0$. There exists an $N \geq 2$ such that for all $n \geq N$, $M_n(\Rnta)$ is contained in an $\epsilon$-neighborhood of 
the spine $\mathcal{S}_t.$
\end{main}

 We establish this result in Section~\ref{sec:LinearSlices}. See Figure~\ref{fig:c=a_n=6_intro}.

\noindent

It nearly goes without saying that all of this work relies on studying the dynamical plane (the Julia sets and polynomial-like maps)--one could consider that Step 0 of the program. In Proposition~\ref{prop:JuliaAnnulus}, we establish an annulus containing the filled Julia set, for any parameters $n\geq 3, c\in\CC, a \in \CC^*$.
Then in Theorem~\ref{thm:babyJulia1}, we show for $n \geq 3, |a| \leq 4, \Arg(a) \neq \pi,$ and $ 0 < c < t=|a|^{1/n}/\max\{4,|a|,|c|\}$, if $n$ is odd then $\Rnca$ restricted to the set of points whose orbits remain in a certain region (a polar rectangle associated with $v_-$)
is topologically conjugate to some quadratic polynomial on its filled Julia set, hence there are homeomorphic copies of quadratic polynomials in the dynamical space. Similarly, in Theorem~\ref{thm:babyJulia2}, we show that for $c < -1$, and for sufficiently large $n$, $\Rnca$ restricted to the set of points whose orbits remain in a certain region (a polar rectangle associated with $v_+$)
is topologically conjugate to some quadratic polynomial on its filled Julia set. We do not complete the remaining steps for either of these two cases but provide this progress as a starting point for future study.

We close the introduction with the organization of the sections. Section \ref{sec:prelim} provides some background information, including the Douady and Hubbard criteria to establish baby $\mandel$'s, and some basic properties of the family $\Rnca$. 
In Section~\ref{sec:dpr} we lay the foundation for parameter spaces by studying the functions $\Rnca$ in the dynamical plane. 
  In Section~\ref{sec:babyMs_Aplane} we prove Main Theorem~\ref{thm:Main_Theorem_APlane}, in the $a$-plane holding $c$ constant. In Section~\ref{sec:LinearSlices}, we study the slices $c=ta$, and prove Main Theorem~\ref{thm:Main_Theorem_C=ta}.

%%%%%%%%%%%%%%%%%%%%%%%%%%%%%%%%%%%%%%%%%%%%%%%%
\section{Preliminaries}
\label{sec:prelim}
%%%%%%%%%%%%%%%%%%%%%%%%%%%%%%%%%%%%%%%%%%%%%%%%

In this section we describe needed tools and results from earlier work.
We begin by explaining how our boundedness locus pictures are drawn.

%-------------------------------------------------------
\subsection{Computer generated images of parameter space in the case of multiple free critical orbits}
To draw a slice of parameter space in the case of two critical orbits, we first fix a relation on the parameters $c,a$ to reduce to one (complex) degree of freedom; for example,  fix a value of $a$ or $c$, or set $c=ta$ for some non-zero $t$. Any one-dimensional (calculable) constraint on the parameters would suffice.
Then, using each of the two critical orbits, we color every point in the picture of the parameter plane in our slice as follows. 

First we assign an (RGB) color to each critical orbit.  Since there are two critical orbits here, $v_+$ and $v_-$, we assign different colors to each value, call them ``color+'' and ``color-'', respectively.  For each parameter value in the picture we test both critical orbits for boundedness and assign a color (and shade) for each orbit.  If the critical orbit appears bounded we assign black. If it escapes, we assign that critical orbit's color, shaded based on rate of escape as is typical; that is, the shade of the color depends on the number of iterations before the orbit escaped a pre-defined escape radius.

Once the testing is complete, each parameter value has two  colors values assigned.  The algorithm then \emph{averages} the two values at each point, resulting in a single assigned (RGB) color for that parameter value.

Therefore a parameter value with both critical orbits bounded will be colored black; a parameter with the critical orbit of $v_+$ bounded while $v_-$ escapes is colored dark color+ and vice-versa is colored dark color-; a parameter with both critical orbits escaping will be colored with the RGB average of the two colors. Note RGB combinations are not the same as paint combinations (for example, red plus green makes yellow), and the colors only truly average if the rate of escapes match. If one escapes more slowly, that color is more intense, so it shades the average toward one of the colors.  Figures \ref{fig:Parameter_planes_Rnca} and \ref{fig:c=a_n=6_intro}, 
%and \ref{fig:Fixed_Crit_Point_intro}) 
give examples of this scheme used to draw an $a$-parameter plane for fixed $c$, and a slice with $c=ta$, respectively.

%-------------------------------------------------------
\subsection{Polynomial-Like Maps}
Douady and Hubbard \cite{DouadyHubbard} show baby $\mandel$'s exist within the parameter planes of other families by showing that a particular family of functions behaves locally like a degree two polynomial.

\begin{defn}\cite{DouadyHubbard} \label{defn:DH1}
A polynomial-like map of degree $d$ is a proper, analytic map of degree $d$, $f:U' \to U$, where $U'$ and $U$ meet the following conditions: both $U'$ and $U$ are open subsets of $\CC$ which are homeomorphic to disks, $U' \subset U$, and $U'$ is relatively compact in $U$.

The \underline{filled Julia set} of a polynomial-like map is the collection of points in $U'$ whose orbits never leave $U'$.
\end{defn}

\begin{figure}
     \includegraphics[scale=.2]{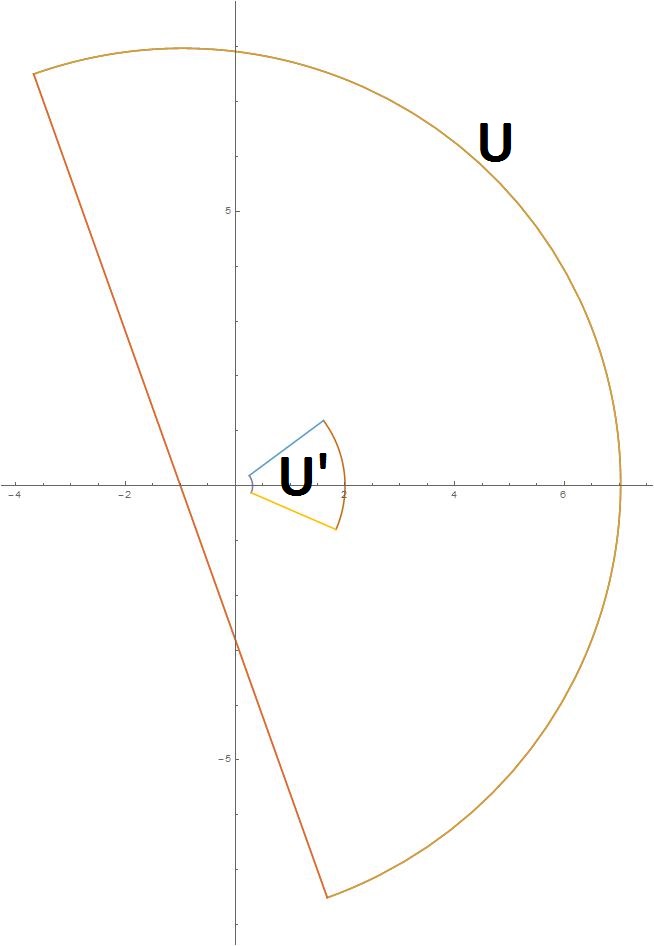}
\caption{\label{fig:UprimeinU} Illustrative example of $U' \subset U$}
    \end{figure}

See Figure~\ref{fig:UprimeinU} for an illustrative example of a $U' \subset U$, with $U'$ a polar rectangle and $U$ a half ellipse, as we will end up having for our cases of interest.

Each polynomial-like map of degree $d$ has exactly $d-1$ critical points. Here we are  only concerned with polynomial-like maps of degree two. A polynomial-like map of degree two is a two-to-one map from $U'$ to $U$, and has a unique critical point in $U'$.  The following shows this is an apt name.

\begin{prop}\cite{DouadyHubbard} \label{prop:DH2}
A polynomial-like map of degree two is topologically conjugate on its filled Julia set to some quadratic polynomial on the polynomial's filled Julia set.
\end{prop}

We call the filled Julia set of the polynomial-like map a ``baby'' Julia set.
Further, as stated by Devaney based on Douady and Hubbard, we have:

\begin{theorem}\cite{DevaneyHalos,DouadyHubbard} \label{thm:DH3}
Suppose we have a family of polynomial-like maps $f_\lambda : U'_\lambda \to U_\lambda$ which satisfy the following:
\begin{enumerate}
    \item The parameter $\lambda$ is contained in an open set in $\CC$ which contains a closed disk $W$;
\item 
 The boundaries of $U'_\lambda$ and $U_\lambda$ both vary analytically as $\lambda$ varies;
\item The map $(\lambda, z) \to f_\lambda(z)$ depends analytically on $\lambda$ and $z$;
\item Each $f_\lambda$ is polynomial-like of degree two and has a unique critical point, $c_\lambda$.
\end{enumerate}
Suppose that for each $\lambda$ in the boundary of $W$ we have that $f_\lambda(c_\lambda) \in U_\lambda-U'_\lambda$, and that $f_\lambda(c_\lambda)$ winds once around $U'_\lambda$ (and therefore once around $c_\lambda$) as $\lambda$ winds once around the boundary of $W$. Then the set of all $\lambda$ for which the orbit of $c_\lambda$ does not escape from $U'_\lambda$ is homeomorphic to the Mandelbrot set.

\end{theorem}

\begin{figure}
  \includegraphics[width=0.45\linewidth,keepaspectratio]{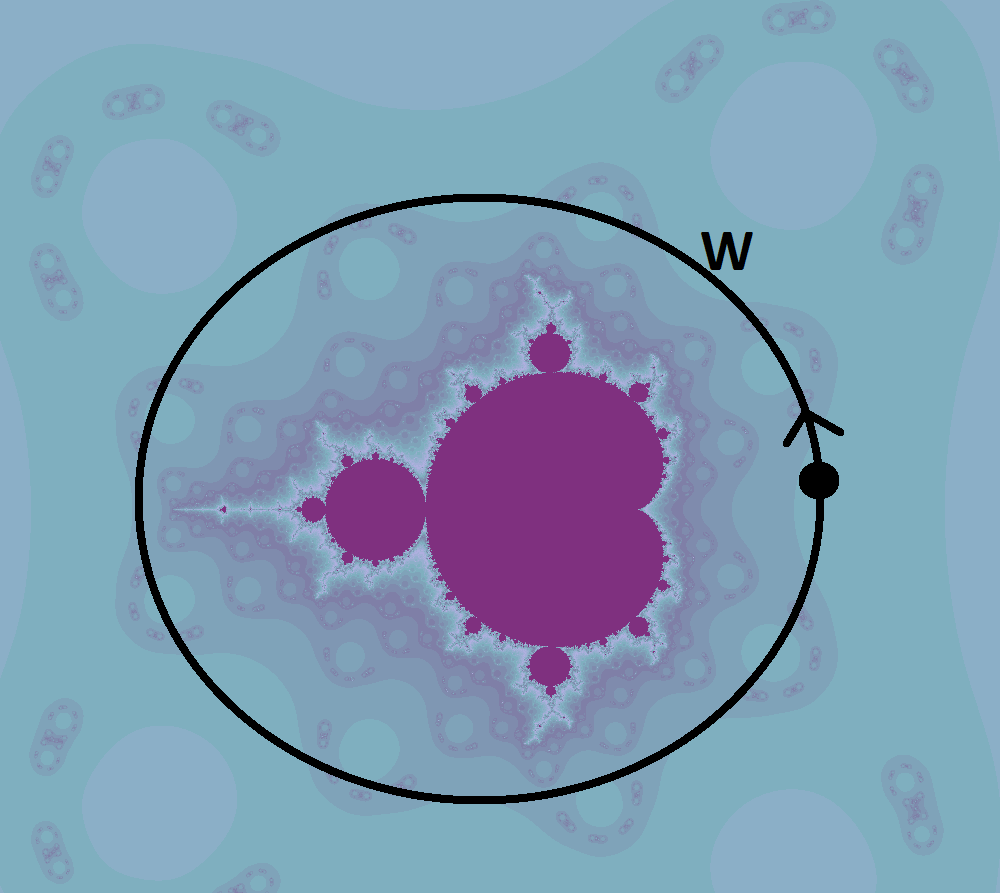} 
  \includegraphics[width=0.45\linewidth,keepaspectratio]{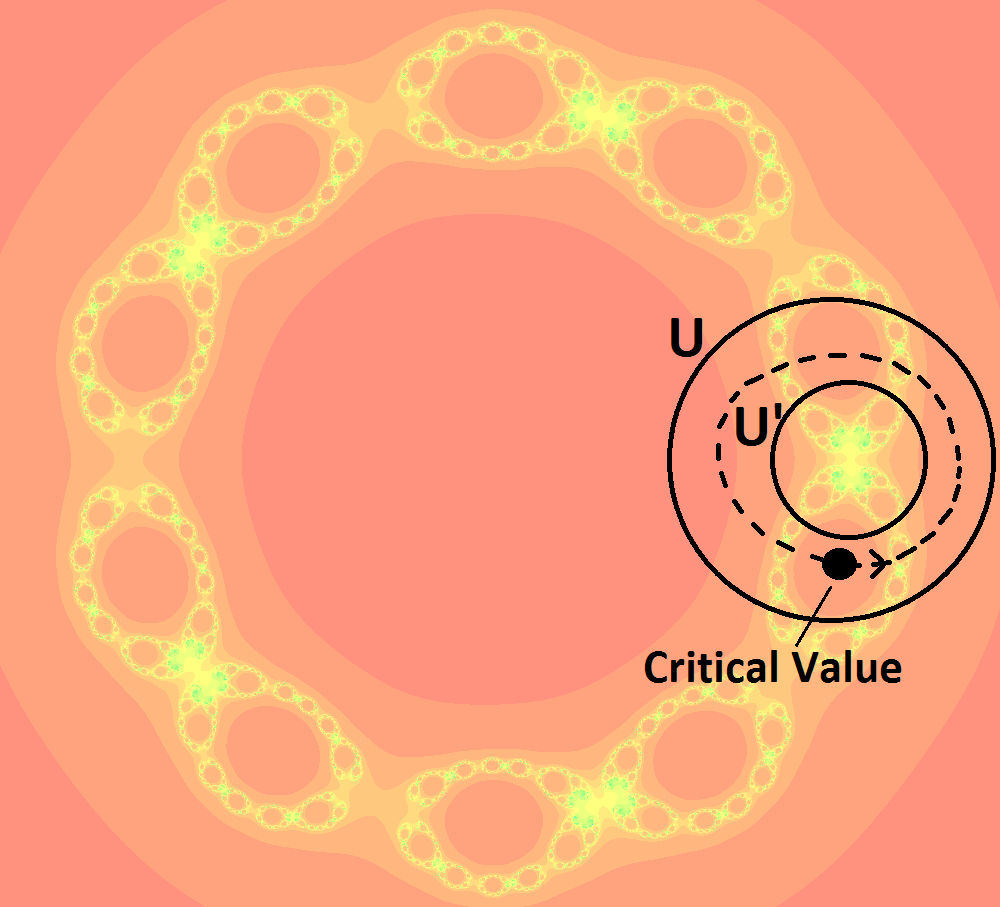}  
\caption{\label{fig:MandelwithW} Establishing the location of a baby Mandelbrot set requires tracking behavior in the dynamical plane as the parameter moves through the parameter plane.}
\end{figure}

Figure~\ref{fig:MandelwithW} illustrates a theoretical loop around a $W$ and corresponding loop around $U - U'$.

This theorem is our method for proving that a given parameter plane contains homeomorphic copies of the Mandelbrot set.  Devaney used this method to show that the $a$-parameter plane of $R_{n,a,0}$ contains $n-1$ baby $\mandel$'s for any $n \geq 3$. In fact, there appear to be many more than just $n-1$ copies, but these are the largest (and with the simplest critical orbit behavior) and so he refers to them as the  ``principal'' copies of the Mandelbrot set.

%-------------------------------------------------------
\subsection{Preliminaries Regarding $\Rnca$}

In \cite{BoydSchulz}, the first author and Schulz provide a number of results on the properties of the family of functions $\Rnca(z)=z^n + \frac{a}{z^n} + c$, as well as on the Julia sets of said family. Two results from this paper are of particular interest here. The first provides information about the location of the filled Julia sets of $\Rnca$ for large $n$.

We use the following notation:

\begin{notation}
$\mathbb{D}(z_0,r) = \left\{z ~ \middle|~ \lvert z - z_0 \rvert < r \right\}$ is the open disc of radius $r$ centered at $z_0$.
\end{notation}

\begin{notation}
$\mathbb{A}(r,R)= \left\{z ~ \middle|~ r < \lvert z \rvert < R \right\}$ is the open annulus centered at the origin lying between the radii $r$ and $R$.
\end{notation}

\begin{cor}\cite{BoydSchulz} \label{cor:BS1}
For any $c \in \CC$ and any $a \in \CC^{*}$, given any $\epsilon > 0$, there is an $N \geq 2$ such that for all $n \geq N$, we have that  the filled Julia set of $\Rnca$ is contained within an annulus near the unit circle: $K(\Rnca) \subset \mathbb{A}(1-\epsilon, 1+\epsilon)$.
\end{cor}

This helps determine a neighborhood of $\Mn(\Rnca)$ for any choice of parameters, for $n$ sufficiently large. This is a fact about the Julia sets, independent of the slice of parameter space (or type of slice) being studied. 

The next result can allow us to locate the centers of potential baby $\mandel$'s in the $a$-plane of $\Rnca$. For $R_{n,a,0}$, Devaney used the fact that one of the baby $\mandel$'s always appears centered on the real axis, along with a rotational argument to locate the $n-1$ principal copies. When we perturb $R_{n,a,0}$ with a nonzero $c$, we need another method in order to locate the baby $\mandel's$. The following lemma is shown in the proof of the last Proposition of \cite{BoydSchulz}, in which the authors find disjoint neighborhoods in the parameter space each containing one fixed critical point.

\begin{lem}\cite{BoydSchulz} \label{cor:BS2}
For any $c \in \CC^*$ and $n \geq 2$, if $|c|\geq 1$ there are $n$ distinct solutions in the family of functions $\Rnca$ to the equations $c \pm 2\sqrt{a} = a^{1/2n}$, and if $|c|<1$ there are $n-1$ distinct solutions.
\end{lem}

 The ``center'' of the Mandelbrot set is $c=0,$ the parameter at which the critical point ($0$) is also a fixed point of the map, so these are natural analogs.

 We will also take advantage of a useful symmetry property of the family of functions $\Rnca$.
\begin{lem}\cite{BoydMitchell}
\label{lem:Involution}
The family of functions $\Rnca$ is symmetric under the involution $h_{a}(z)=\frac{a^{1/n}}{z}$.
\end{lem}

%%%%%%%%%%%%%%%%%%%%%%%%%%%%%%%%%%%%%%%%%%%%%%%%
\section{General Dynamical Plane Results}
\label{sec:dpr}
%%%%%%%%%%%%%%%%%%%%%%%%%%%%%%%%%%%%%%%%%%%%%%%%

In this section, we establish a few properties concerning the dynamics of $\Rnca$ and the sets $U$ and $U'$ used to establish polynomial-like behavior, which will remain useful throughout the rest of this paper and in future study. For flexibility, note we make no default assumptions on the parameters, and only when necessary bring in minimal restrictions, and that is simply $n \geq 3$ and $c,a$ complex with $a\neq 0$, and recall $\psi=\Arg(a)$. 

%====================================================
\subsection{The sets $U_{a,k}$'s and their images for any $\Rnca$.}
%====================================================

We have pointed out that this family has $2n$ critical points and two critical values, $v_{\pm} = c \pm 2\sqrt{a}$.
In addition, for this family, $R$ maps half the critical points (counting multiplicity) to $v_+$ and the other half to $v_-$. Specifically, a simple calculation yields:
\begin{lem} The critical points of $R_{n,a,c}$ satisfy
\label{lem:critptsmapping}
\begin{equation}
\label{eqn:allcritpts}
a^{1/2n} = |a|^{1/2n}e^{i(\frac{\psi + 2k\pi}{2n})} 
\;\;\; \text{for} \;\;\; k=0,1,\ldots,2n-1,
\end{equation}
where $\psi= \Arg(a)$. Moreover, the $n$ critical points defined by $k$ even map to $v_+=c+2\sqrt{a}$ while the $n$ critical points with $k$ odd map to $v_-=c-2\sqrt{a}$.
\end{lem}
That is, $\crpt_{2j-1} = |a|^{1/2n}e^{i(\frac{\psi + 2(2j-1)\pi}{2n})}$
for $j=1,2,3,...,n,$ map to $v_-$ and 
$\crpt_{2j}=|a|^{1/2n}e^{i(\frac{\psi + 2(2j)\pi}{2n})} $ for
$j=0,1,2,...,n-1,$ map to $v_+.$ 

Note also that $\crpt_0 = |a|^{1/2n} e^{i\frac{\psi}{2n}}$ and 
 $\crpt_n = -|a|^{1/2n} e^{i\frac{\psi}{2n}}$.

\begin{proof} For $k\in \{0, 1,\ldots, 2n-1\}$, let
$\crpt_k = |a|^{1/2n} e^{i(\psi+2k\pi)/2n}$.
We want $R(\crpt_k)= c\pm 2\sqrt{a}$ so will examine whether
$(R(\crpt_k) -c)/(\sqrt{a}) = \pm 2$.

Now, 
\begin{align*}
& \frac{R( |a|^{1/2n}e^{i(\frac{\psi + 2(k)\pi}{2n})}) - c}{\sqrt{a}} 
\\ & 
= \frac{ (|a|^{1/2n}e^{i(\frac{\psi + 2(k)\pi}{2n})})^n }{|a|^{1/2}e^{i\psi/2}}
+ \frac{a}{(|a|^{1/2n}e^{i(\frac{\psi + 2(k)\pi}{2n})})^n
|a|^{1/2}e^{i\psi/2} }
\\ & 
= \frac{ (|a|^{1/2}e^{i(\frac{\psi + 2(k)\pi}{2})}) }{|a|^{1/2}e^{i\psi/2}}
+ \frac{|a|e^{i\psi}}{(|a|^{1/2}e^{i(\frac{\psi + 2(k)\pi}{2})})
|a|^{1/2}e^{i\psi/2} }
\\ & 
= \frac{ (e^{i(\frac{2(k)\pi}{2})}) }{1}
+ \frac{1 }{(e^{i(\frac{2(k)\pi}{2})})  }
\\ & 
= e^{ik\pi} + e^{-ik\pi},
\end{align*}
which is equal to $2$ when $k$ is even and $-2$ when $k$ is odd.
\end{proof}

Now, to apply Theorem~\ref{thm:DH3} to prove that we have a baby $\mandel$, we will need to create the two sets, $U'$ and $U$, and prove that they satisfy several conditions. 
We will definitions are similar to those from \cite{DevaneyHalos} and \cite{BoydMitchell}.

\begin{defn}
\label{defn:U_ak_prime}
For the family of functions $\Rnca$, for any $a$ with $|a| < 4^n$, let $\psi = \Arg(a)$ and for each $k=0,\ldots,2n-1,$ define 
$U_{a,k}'$ (which also depends on $n$),  
as the polar rectangle:
 $$U'_{a,k} = \left\lbrace z \;\middle|\; \frac{|a|^{\frac{1}{n}}}{2} < |z| < 2, \; 
\left| \Arg(z) - \left( \frac{\psi+2k\pi}{2n} \right) \right| < \frac{\pi}{2n} 
\right\rbrace.$$
Also, set $U_a^+$ to be the image of $U'_{a,0},$ and $U_a^-$ to be the image of $U'_{a,1}$:
$$U_{a}^+ = \Rnca (U'_{a,0}), \ \text{ and } \  U_{a}^- = \Rnca (U'_{a,1}).$$
\end{defn}

\begin{figure}
     \includegraphics[scale=0.7]{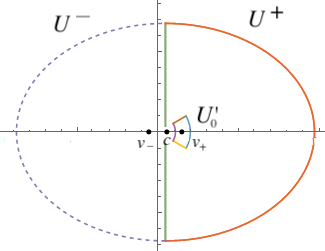}
\caption{\label{fig:GoodEllipseLabelled} Diagram example of sets $U'_0, U^+, U^-$ for parameters $c>0,a<0$ such that in this case, $U'_0 \subset U^+$.}
    \end{figure}

So $U'_{a,k}$ is the polar rectangle with arguments centered at $(\psi + 2k\pi)/2n$ and then plus or minus $\pi/2n$.
Note $U'_{a,k}$ does not depend on $c$, but of course its image does. 
Figure~\ref{fig:GoodEllipseLabelled} shows an example of sets $U'_0, U^+, U^-$ for real parameters $c,a$ such that in this case, $U'_0 \subset U^+.$

\begin{figure}
\centering
\begin{subfigure}{.5\textwidth}
  \centering
\includegraphics[width=1.0\textwidth,keepaspectratio]{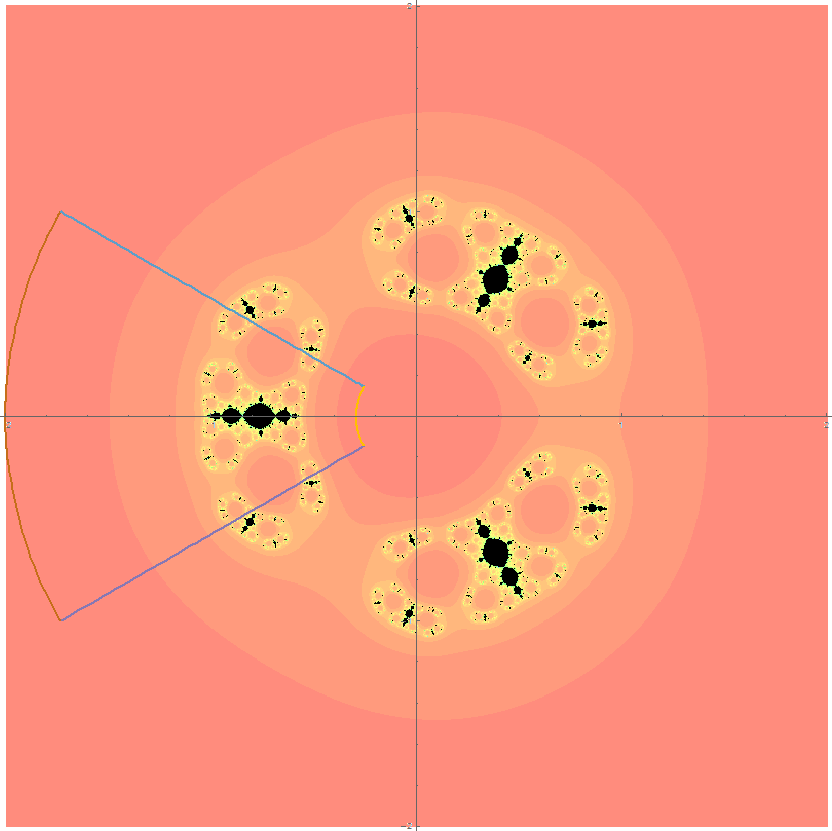}
  \caption{$U'_{a,3}$ overlaid on the Julia set of $R_{3,0.2,0.25}$}
  \label{fig:JuliaOverlay1} 
\end{subfigure}%
\begin{subfigure}{.5\textwidth}
  \centering
\includegraphics[width=1.0\textwidth,keepaspectratio]{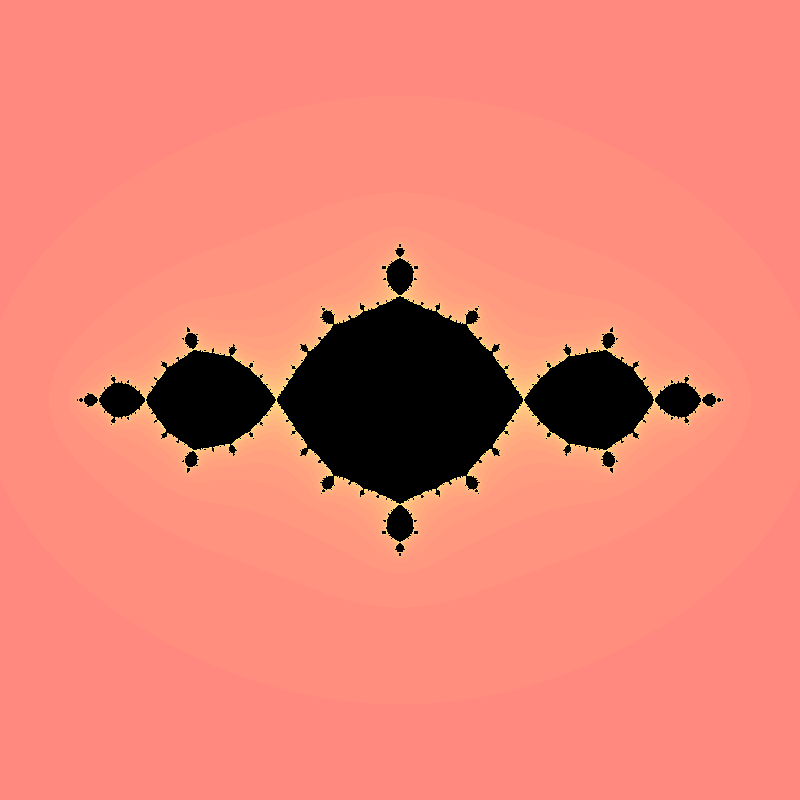}
  \caption{The filled Julia set of $P_{-1}(z)=z^2-1$}
  \label{fig:Basilica} 
\end{subfigure}
\caption{Left is an example of a $U'_a$ containing a baby Julia set, right is the filled Julia set of $z^2-1$}
\label{fig:OverlayContrast1}
\end{figure}

Figure \ref{fig:JuliaOverlay1} shows an example of a filled Julia set of $\Rnca$ with an overlay of  $U'_{a,3}$. Figure \ref{fig:Basilica} shows the filled Julia set  $P_{-1}(z)=z^2-1$. One sees why these subsets are called ``baby'' quadratic Julia sets.

%---------------------------------------------------
\begin{defn}
\label{defn:ellipseEE}
 Let $\EE$ denote the ellipse parameterized by
$$
x(\theta) = \left( 2^n + {|a|}/{2}\right) \cos(\theta), \ \ \text{ and }
y(\theta) = \left( 2^n - {|a|}/{2}\right) \sin(\theta),
$$
then shifted so that it is centered at $c$ and rotated counter-clockwise by ${\psi}/{2}$.
\end{defn}

Note before the shift and rotation, the ellipse $\EE$ is centered at $0$,  its semi-major axis lies along the $x$-axis and has length $\left( 2^n + {|a|}/{2}\right)$, and its semi-minor axis lies along the $y$-axis and has length $\left( 2^n - {|a|}/{2}\right)$. 

By Lemma 6 of \cite{BoydMitchell}, the foci of $\EE$ are at the critical values $v_{\pm}.$

\begin{prop} \label{prop:Uellipse}
For each $k$ even in $\{ 0,1, \ldots, 2n-1\}$, 
% $U_a = \Rnca(U'_{a})$ 
$U_a^+ = \Rnca(U'_{a,k})$ 
and this set is one half of the ellipse $\EE$, including the minor axis and the critical value $v_+$. 
Moreover, $\Rnca$ maps 
$U'_{a,k}$ $2:1$ onto $U_a^+$ for each even $k$.

Similarly, 
$U_a^- = \Rnca(U'_{a,k})$ for each odd $k$
is the other half of the ellipse $\EE$, including the minor axis and the other critical value $v_-$, and $\Rnca$ maps 
 $U'_{a,k}$ two-to-one onto $U_a^-$ for each odd $k$.
\end{prop}

The $k=0$ case is established in \cite{BoydMitchell},  in the first section of results, in proofs of the Lemmas ``$U$ is half an ellipse centered at $c$ and rotated by $\psi/2$.'' and ``$v_{\pm}$ are the foci of $\EE$''. The authors stated some assumptions on the parameters, but a careful review reveals those were not used in the proofs of the relevant lemmas on $\EE$. 

\begin{proof}
First we examine the images of the outer and inner boundary arcs of any
 $U'_{a,k}$. 
 Note that, due to the involution symmetry, both the outer arc which lies on $|z|=2$ and the inner arc which lies on $|z|=\frac{a^{1/n}}{2}$ are mapped  by $\Rnca$ to the same curve. 

Consider $\Rnca (2e^{i\theta})$ as $\theta$ goes from $0$ to $2\pi$. Since we are considering all values of $\theta$, we can add a rotation to the angle of $\theta$ and still end up with the same image. Now, a calculation (found in \cite{BoydMitchell}), yields:
\begin{equation*}
\label{eqn:RofEE}
 \Rnca \left( 2e^{i(\theta + \frac{\psi}{2n})}\right) = 
e^{i\frac{\psi}{2}} \left(x(n\theta) + iy(n\theta) \right) + c,
\end{equation*}
where $x$ and $y$ are the same as in Definition~\ref{defn:ellipseEE}. Thus 
$\Rnca (2e^{i\theta})$ as $\theta \colon 0 \rightarrow 2\pi$ is exactly $\EE$ traversed $n$ times.

Now, the image under $\Rnca$ of both circles of radius $2$ and $\frac{a^{1/n}}{2}$ covers $\EE$ exactly $n$ times, but the boundary of  
any $U'_{a,k}$
only occupies $\frac{1}{2n}^{th}$ of each circle. Thus the image of the inner and outer arcs of 
each $U'_{a,k}$
is exactly half of $\EE$.

Next we evaluate the images under $\Rnca$ of the boundary rays  of
 $U'_{a,k}$: 
$R_{n,c,a}\left( r*\exp \left(i\frac{\psi + 2k\pi \pm \pi}{2n} \right) \right) = $
\begin{eqnarray*}
& = & \left( r \cdot \exp \left(i\frac{\psi + 2k\pi  \pm \pi}{2n} \right) \right) ^{n}+\frac{a}{\left( r\cdot \exp \left(i\frac{\psi + 2k\pi  \pm \pi}{2n} \right) \right) ^{n}}+c\\
& = & r^{n}\cdot\exp \left(i\frac{\psi + 2k\pi  \pm \pi}{2} \right)+\frac{|a|e^{i\psi}}{r^{n}\cdot\exp \left(i\frac{\psi + 2k\pi  \pm \pi}{2} \right)}+c\\
& = & r^{n}\cdot\exp \left(i\frac{\psi + 2k\pi  \pm \pi}{2} \right)+\frac{|a|}{r^{n}}\cdot\exp \left(i\frac{\psi - 2k\pi  \mp \pi}{2} \right)+c\\
\end{eqnarray*}
\begin{eqnarray*}
& = &
e^{i\frac{\psi}{2}} 
\left(r^{n}\cdot\exp \left( i\frac{2k\pi  \pm \pi}{2} \right
) - \frac{|a|}{r^{n}}\cdot\exp \left( i\frac{2k\pi  \pm \pi}{2} \right) \right)+c \\
& = &
e^{i\frac{\psi}{2}} 
 \exp \left(i \frac{2k \pi }{2} \right) 
\exp \left(i \frac{\pm \pi}{2} \right)
\left(r^{n} - \frac{|a|}{r^{n}}\right)i+c\\
& = &
 e^{i k \pi }  e^{i\frac{\psi}{2}} \left(r^{n} - \frac{|a|}{r^{n}}\right)i+c,
\end{eqnarray*}
and $e^{i k \pi }=(-1)^k$.

Each of these is a line segment along the imaginary axis which is then rotated by $\psi/{2}$ and shifted by $c$. The images of the two boundary rays correspond to the same line segment, traversed in opposite directions as $r: \frac{a^{1/n}}{2} \mapsto 2$. Before rotation by ${\psi}/{2}$ and translation by $c$, the line segment $\left( r^n - \frac{|a|}{r^n} \right)i$ as $r: \frac{a^{1/n}}{2} \mapsto 2$ is the line segment between $\left( 2^n - \frac{|a|}{2^n} \right)i$ and $-\left( 2^n - \frac{|a|}{2^n} \right)i$. After rotation and translation, this is the minor axis of  $\EE$. 

 Recall from Lemma~\ref{lem:critptsmapping}, the $k$ even critical points, those located in the $U_{a,k}'$s for $k$ even map to $v_+$ in $U_a^+$, and the $k$ odd ones map to $v_-$ in $U_a-$.
 
 Thus, each of
 $U_a^{\pm}$
 is a semi-ellipse, divided along the minor axis, centered at $c$,  rotated by ${\psi}/{2}$, and including as focus $v_{\pm}$.
 Moreover, $\Rnca$ maps 
 $U'_{a,k}$ 
two-to-one  onto $U_a^+$ for $k$ even, and two-to-one onto $U_a^-$ for $k$ odd.
\end{proof}

%%%%%%%%%%%%%%%%%%%%%%%%%%%%%%%%%%%
\subsection{The Escape Radius for $\Rnca$.}
Now we  establish some properties of the filled Julia sets of $\Rnca$ which are true with almost no restrictions on the parameters. These are generalizations of results in \cite{BoydSchulz}. To that end we only assume $a \in \CC^*$, $c \in \CC$, and $n \geq 3$. Throughout this paper we ignore the cases $n=1$ and $n=2$, as the dynamics are less interesting and not in line with the results for any larger values of $n$, respectively (see \cite{DevaneyCase2}.)

For $n\geq 2$, the point at infinity is a super-attracting fixed point, meaning that the derivative when evaluated at infinity is zero (as usual, view the map as a map of the Riemann sphere and calculate the derivative at infinity by conjugating with the function $z \to {1}/{z}$ and calculating the derivative at zero). 
%This fact does not hold true when $n=1$ and is one reason that case is left out.
 Infinity being super-attracting implies that there is an escape radius beyond which the orbits of $z$ must tend towards infinity. Therefore each Julia set is bounded, and avoids some neighborhood of infinity. In this first result we calculate a neighborhood of the origin, as a function of the parameters, in which the filled Julia set must lie. This is an adaptation of a lemma from \cite{BoydSchulz}. We note in \cite{BoydMitchell}, in both studies on the $a$ plane and $c$ plane parameter restrictions are utilized so that the escape radius is always $2$. 

\begin{lem} \label{lem:OuterEscape}
For any $c \in \CC$, $a \in \CC^*$, and integer $n \geq 3$, set $s=\max\{4,|c|,|a|\}$. Then for any $|z| \geq s$, the orbit of $z$ under $\Rnca$ escapes to infinity.
\end{lem}

\begin{proof}
Fix $c \in \CC$, $a \in \CC$, and $n \geq 3$. Let $|z| \geq s$. Note $s^{n-2} \geq 4$ for any choice of $n \geq 3$. We show by induction $|\Rnca^{m}(z)| > s^{m}$ for all $m \geq 1$.

First, 
\begin{align*}
& |\Rnca(z)|  = \left | z^n + c + \frac{a}{z^n}\right | 
 \geq |z|^n - |c| - \frac{|a|}{|z|^n} 
 > 3s - s - \frac{|a|}{3s} 
 \\& 
 > 3s - s - \frac{|a|}{s}
 = 2s - \frac{|a|}{s}
 = s + \left (s - \frac{|a|}{s} \right ) 
 > s.
\end{align*}
Now, suppose for some $m \geq 1$, $|\Rnca^m(z)| > s^m$.  Then,
\begin{align*}
& |\Rnca^{m+1}(z)|  
= \left | (\Rnca^m(z))^n + c + \frac{a}{(\Rnca^m(z))^n} \right |
\\  & 
\geq   | \Rnca^m(z)|^n - |c| - \frac{|a|}{ | \Rnca^m(z)|^n}
> s^{mn} - s - \frac{|a|}{s^{mn}}
\\ & 
\geq s^{mn-m+1} - s - \frac{|a|}{s^{mn}}
= s^{m(n-1)+1} - s - \frac{|a|}{s^{mn}}
\\ &
 = s^{m+1} \left ( s^{m(n-2)} - s^{-m} - \frac{|a|}{s^{mn+m+1}} \right )
 > s^{m+1} (4^{m}-1-1) 
 > s^{m+1}.
\end{align*}
Thus by induction, $|\Rnca^m(z)| > s^m$ for all $m \geq 1$. Since $s \geq 4$, 
the orbit of $z$ under $\Rnca$ escapes to infinity. Thus we have that for any choice of $z$ such that $|z|>s$, the orbit of $z$ under $\Rnca$ will escape to infinity.
\end{proof}

Next we notice that, since $0$ is a preimage of infinity, there must be a neighborhood of $0$ which maps beyond the escape radius established in Lemma \ref{lem:OuterEscape}. Thus there must be a neighborhood of $0$ which the filled Julia set avoids as well, which we find using the involution from the previous section. Again, this is a more general result than those in prior publications.

\begin{lem} \label{lem:InnerEscape}
For any $c \in \CC$, any $a \in \CC^*$, and any integer $n \geq 3$, set $t=\dfrac{|a|^{1/n}}{s}$, where $s$ is defined as in Lemma \ref{lem:OuterEscape}. Then $t<s$, and for any $|z| \leq t$ the orbit of $z$ under $\Rnca$ escapes to infinity.
\end{lem}
\begin{proof}
First, we show $t<s$, and $t<1$ by considering three cases. Case 1: suppose $|a|<1$. Then $|a|^{1/n}<1$, and since $s\geq 4$ we have $t<1$, and thus $t<s$. Case 2: suppose $|a|=1$. Then $|a|^{1/n} = 1$ and Lemma \ref{lem:OuterEscape}, $s\geq 1$. Hence $t=1/s$, so $t<s$ and $t<1$. Case 3: suppose $|a|> 1$. Then $|a|^{1/n}< |a|$, so $t= < |a|/s$. Since $s\geq |a|$, we have $t<1$ again, so $t<s$. Thus $t<s$ and $t<1$.

Now, let $z$ be such that $|z| \leq t$. Then $$|H_a(z)|=\frac{|a|^{1/n}}{|z|} \geq \frac{|a|^{1/n}}{t}
=s.$$ 
Thus, by Lemma \ref{lem:OuterEscape}, the orbit of $H_a(z)$ under $\Rnca$ escapes to infinity. Now, the involution has the property that $\Rnca(z)=\Rnca(H_a(z))$, so since the orbit of $H_a(z)$ escapes to infinity, so does the orbit of $z$.
\end{proof}

Combining Lemmas \ref{lem:OuterEscape} and \ref{lem:InnerEscape} gives us an annulus containing the filled Julia set $K(\Rnca)$, for the full family of interest:

\begin{prop} \label{prop:JuliaAnnulus}
Let $c \in \CC$, $a \in \CC^*$, and $n \geq 3$. Let $s=\max\{4,|c|,|a|\}$ and $t=\dfrac{|a|^{1/n}}{s}$. Then 
 $K(\Rnca) \subset \mathbb{A}(t,s)$.
\end{prop}

%%%%%%%%%%%%%%%%%%%%%%%%%%%%%%%%%%%%%%%%%%%%%%%%
\section{Baby Julia Sets and Baby Mandelbrot Sets in the case of holding $c$ constant}
\label{sec:babyMs_Aplane}
%%%%%%%%%%%%%%%%%%%%%%%%%%%%%%%%%%%%%%%%%%%%%%%%

In this section we present our first example of following the program that we have set forth. In the case that the parameter space is the $a$-plane with the parameters $c$ and $n$ held constant, we will investigate the presence of baby Julia sets in dynamical plane, the boundedness locus $\Mn(\Rnca)$, and the presence of baby $\mandel$'s, ultimately establishing Theorem~\ref{thm:Main_Theorem_APlane}.

%%%%%%%%%%%%%%%%%%%%%%%%%%%%%%%%%%

We have established definitions and descriptions for $U'_{a,k}$ (a polar rectangle) and $U_a^{\pm}$.
We know each $U'_{a,k}$ maps two-to-one onto $U_a^{\pm}$, but we have not established that $U'_{a,k} \subset U_a^{\pm}$ for any values of $n$, $a$ or $k$. Not every map in this family has a subset on which it is polynomial-like. 
In this section we show that $U'_{a,k} \subset  U_a^{\pm}$ for a few different types of parameter restrictions, establishing the existence of some baby Julia sets in the Julia sets of $\Rnca$. 

We also find a polar rectangle in the $a$-plane that contains $\Mn(\Rnca)$ for a  wide range of parameter values.

This section is split into cases based on the value of the parameter $c$, as the type and quantity of principal baby $\mandel$'s that appear in the $a$-plane of $\Rnca$ changes depending on whether $c$ is very small or relatively large. This change in behavior happens around $|c|=1$. As we saw in Section~\ref{sec:prelim}, when $|c| \geq 1$ there are $n$ solutions to the equations $c \pm 2\sqrt{a} = a^{1/2n}$, which correspond to $n$ principal baby $\mandel$'s. However, when $|c| < 1$, some of these principal baby $\mandel$'s start to merge and deform, resulting in shapes which have some similarities to $\mandel$, but which are noticeably different.

Primarily we focus on the dynamics when $|c| > 1$, finding baby $\mandel$'s there. 
The case of $c\in[-1,1]$ was considered in \cite{BoydMitchell}.  Nevertheless, in the first subsection below we do establish some polynomial-like behavior for some smaller ranges of $c$, which are distinct from the results of \cite{BoydMitchell}.

%%%%%%%%%%%%%%%%%%%%%%%%%%%%%%%%%%
\subsection{Dynamical Plane Results for Small $c$}

To show that some $U'_{a,k}$ is contained in $U_a^{\pm}$ given certain restrictions,  we first we show  each $U'_{a,k} \subset \EE$, and then show under some additional restrictions that a $U'_{a,k}$ is contained in its image, one of $U_a^{\pm}$, by showing $U'_{a,k}$ does not touch the minor axis of $\EE$ (and $k$ is such that $U'_{a,k}$ intersects its image, i.e., it's in the right half of $\EE$).

\begin{prop} \label{prop:SmallC}
Let $|c|<t$ as defined in Lemma \ref{lem:InnerEscape} and let $|a|\leq 4$. Then for any $n \geq 3$, $\DD(0,2) \subset \EE$ and for all $k\in \{ 0, \ldots, 2n-1\},$  $U'_{a,k} \subset \mathcal{E}$.
\end{prop}

\begin{proof}
We have $|c|<t=\dfrac{|a|^{1/n}}{s}$, where $s=\max\{4,|c|,|a|\}$, and $|a|\leq 4$. Then
$|c|<\dfrac{|a|^{1/n}}{4}$. 
The shortest distance from the center, $c$, of $\mathcal{E}$ to the boundary is given by the length of the semi-minor axis, which is $2^n-\dfrac{|a|}{2^n}$. Since $n\geq 3$, $2^{n+1}-\dfrac{1}{2^{n-1}} > 2^n-\dfrac{1}{2^{n-2}}$ for any choice of $n$.\\
Thus,
$$2^n-\frac{|a|}{2^n} \geq 2^n - \frac{4}{2^n} = 2^n-\frac{1}{2^{n-2}} \geq 8-\frac{1}{2} = \frac{15}{2}.$$
Next, suppose $|z|=2$. Then $|c-z| \leq |c|+|z| < \dfrac{|a|^{1/n}}{s} + 2 < 1+2=3$. We now have that the distance from from $c$ to any point on the circle $\{z \mid |z|=2 \}$ is less than $3$, and the length of the semi-minor axis of $\mathcal{E}$ is greater than $\frac{15}{2}$. Thus $\mathcal{E}$ contains $\{z \mid |z| \leq 2 \}$, and therefore contains every $U'_{a,k}$.
\end{proof}

Now if we further restrict to the case of $c$ real and positive, we can show $U'_{a,n} \subset U_a^-$ for $n\geq 3$ an odd number. In this case, with $n$ odd recall we have $R(U'_{a,n})=U_a^-$.  (The motivated reader might consider attempting to generalize the reminder of this subsection, that is, examine the $n$ even case for $c$ real, and then add an argument rotation to $c$ and examine that rotation to determine which $U'_{a,k}$ contains the bounded critical value, and then try to establish polynomial-like behavior for those complex $c$ values).  

First, a lemma to help us envision the relative positions of the key regions and values in this case.

\begin{lem}
\label{lem:sign_of_vminus}
Let $n\geq 3$, and $c, a \in \RR^+$ with $0<a\leq 4$  and $0<c<t=\dfrac{a^{1/n}}{s}$. Then $v_+>0$ and 
\begin{enumerate}
    \item  $v_- >0$ for all $a$ with $(1/4)^\frac{n}{n-1}< a < 4,$ any $c$ as restricted above.
\item $v_- >0 $ for $0< a < (1/4)^\frac{n}{n-1}$ if $0< c < 2\sqrt{a}$ (and such $c$ occur, as $0<a$ so $0< 2\sqrt{a}$).
\item $v_- \leq 0$ for $0< a < (1/4)^\frac{n}{n-1}$ if $ 2\sqrt{a} \leq c < a^{1/n}/4$ (and such $c$ also occur).
\end{enumerate}
\end{lem}

\begin{proof}
    %--------
 For $c>0$ and $a>0$, we know $v_+ = c + 2\sqrt{a}>0$.  Now we consider, for these parameter restrictions, the sign of $v_-$, and 
 $v_-<0$ if and only if $c < 2\sqrt{a} $.

 We have assumed that $|a| \leq 4$, so certainly we have $|a|^{1/2n} < 4$ 
 (because if $|a| \in (1,4]$ then its roots are smaller, and if $|a| \geq 1$ then roots are also bounded by $1$, hence by $4$). 
  From this we can say that $\frac{|a|^{1/n}}{4} < |a|^{1/2n}$ and therefore $c < |a|^{1/2n}$ for all available choices of $a$ (also note then  $c<1$).
 
 \textit{Case i}: Suppose $1 \leq a<4$. Then $1 \leq \sqrt{a} < 2,$ so $2 \leq 2\sqrt{a} < 4$, and since $c<1<2$ we have  $c < 2\sqrt{a}$, so for $a\in [1,4)$ we have $v_- <0$.  
 
\textit{Case ii}:  Suppose $0<a < 1$. Then $0 < a < \sqrt{a}  < a^{1/2n} < 1$ so $0< 2a < 2\sqrt{a} < 2a^{1/2n} < 2$, and  
$ c < a^{1/2n}$ means $ 2c < 2a^{1/2n}$.
Then 
note, we can have $0 \leq v_-$ i.e., $2\sqrt{a} \leq c$ if $2\sqrt{a} \leq c < 2c < 2a^{1/2n}$ or 
$\sqrt{a} \leq c/2 < c < a^{1/2n}$ or $a \leq c^2/4 < c^2 < a^{1/n}$. We know $c< a^{1/n} / 4$ but $c<1$ also 
so $c^2 < c< a^{1/n} / 4$, 
hence we are considering 
\begin{equation}
\label{eqn:a_and_c_bounds}
a \leq c^2/4 < c^2 < c <a^{1/n} / 4 < a^{1/n}.
%(*) 
\end{equation}
Can we have $a < a^{1/n}/4$? Sure if $a^n < a / 4^n$ or $a^n - a/4^n < 0,$ or $a (a^{n-1} - 1/4^n) < 0$, and since $a>0$ we need $(a^{n-1} - 1/4^n) < 0$ or $a^{n-1} < 1/4^n =(1/4)^n $ so $a < (1/4)^\frac{n}{n-1}$. 
In summary, for $0< a < (1/4)^\frac{n}{n-1}$ then
$a < a^{1/n}/4$  and so there is room to allow a range of $c$ for which inequality~\ref{eqn:a_and_c_bounds} can occur so for $0< a < (1/4)^\frac{n}{n-1}$ and $c$ such that $c \in [2\sqrt{a}, a^{1/n}/4)$ (which is not an empty interval) we have $v_- \geq 0$, with $v_-=0$ for $c=2\sqrt{a}.$ Then for $(1/4)^\frac{n}{n-1}< a < 1,$ we have $a^{1/n}/4 < a$ so $v_->0$. 

Hence we have established (1) by combining Case i and part of Case ii, and (2) and (3) from Case ii. 
\end{proof}

\begin{prop} \label{prop:FirstResult}
Let $n\geq 3$,  $c \in \RR^+$, $0<c<t=\dfrac{|a|^{1/n}}{s}$, $0<|a|\leq 4$,  and $\Arg(a)=\psi \in ( -\pi, \pi]$. Then for $n$ odd, $U'_{a,n} \subset U_a^-$.  
\end{prop}

\begin{proof}
First, these restrictions on $c$, $a$, and $n$ satisfy the requirements of Proposition \ref{prop:SmallC}, so $U'_{a,k} \subset \mathcal{E}$ for every $k$. Hence $U'_{a,n} \subset \EE$ (and $\Rnca$ maps $U'_{a,n}$ $2:1$ onto its image $U_a^-$, by Proposition~\ref{prop:Uellipse}).

Since $c \in \RR^+$, for the moment let's consider an $a$ so that $a \in \RR^+$ as well, so $a=|a|$ and $\psi = \Arg(a)=0$. Then we get an $\EE$ 
%which is centered on the real axis 
with its major axis lying on $\RR$ and its minor axis, which makes up part of the boundary of $U_a^{-}$, running vertically through $c>0$. 

Since $a>0$, we also have $\sqrt{a}$ is real so $\pm 2\sqrt{a}$ are two points on the real axis, $\sqrt{a}>0$ and  $-\sqrt{a}<0$. Hence the foci of $\EE$, the critical values $v_{\pm}$, lie on the real axis and from the above lemma we know $v_+>0$, and the sign of $v_-$ can vary with $v_->0$ for larger $a$, or for $c$ small relative to $a$, but $v_-$ can be positive if $a$ is small and $c$ is a bit larger relative to $a$. 
 Since $a>0$, $v_-=c-2\sqrt{a}< c$, so $v_-$ lies in the semi-ellipse half of $\EE$ intersecting the left half plane of $\CC$, $\{z : \Re(z) < 0\}$, and note $v_+$ lies in the right half of $\EE$ lying completely in the right half plane $\{z : \Re(z) >0\}.$

 Now, from the definition we see that $U'_{a,n}$ is angularly centered around the (negative) real axis, with the two ray boundaries at arguments $\pi \pm \frac{\pi}{2n}$. Since $U'_{a,n}$ lies in the left half plane, is a subset of $\EE$, and $c>0$ is the center of $\EE,$ we have $U'_{a,n}$ lying entirely in $U_a^-$ and disjoint from $U_a^+$, again assuming $a>0$.

Next we let $\psi=\Arg(a)$ increase from $0$ to $\pi$. In doing this $\mathcal{E}$ is rotated by ${\psi}/{2}$, so as $\psi$ increases, the argument of each point in $U_a^{\-}$ with respect to $c$ increases by ${\psi}/{2}$. On the other hand, as $\psi$ increases the argument of each point in any $U'_{a,k}$ increases by only $\psi/2n$ so at most ${\pi}/{6}$ since $n \geq 3$. 
Thus, the ray boundary $\pi + \pi/2n$ of $U'_{a,n}$ can't catch up with the minor axis, but we have to worry about whether the minor axis catches up with the ray boundary $\pi - \pi/2n$. 
This can happen when $\psi=\pi$, so $U_a^{-}$ has been rotated by $\pi/2,$ so that the minor axis of $\EE$ is now in the real line (still centered at $c>0$), and $\pi - \pi/(2n) + \psi = \pi - \pi/2n + \pi/(2n)=\pi$.

Since the argument of the $U'$s and the $U$s both increase linearly, this is the first and only time that the boundary of $U_a^-$ touches the boundary of $U'_{a,n}$. Thus for $0\leq \psi <\pi$ we have $U'_{a,n} \subset U_a^-$. Furthermore, when $\psi=\pi$ the boundaries of $U'_{a,n}$ and $U_a^-$ touch, but both $U'_{a,n}$ and $U_a^-$ are open so we still have $U'\subset U$ (and we still have a proper subset due to the other three boundaries of the $U'$s).

Finally we let $\psi$ go from $0$ to $-\pi$. This works the same as the previous case, with the only difference being a clockwise rotation instead of counter clockwise. When $\psi=-\pi$, $U_a^-$ has been rotated by $-{\pi}/{2}$, so the minor axis of $\mathcal{E}$ lies on the real axis and the $\pi + \pi/2n$ ray boundary of $U'_{a,n}$ has argument $\pi$. 

Again, the argument of each $U'$ and $U$ both change linearly so this is the first and only time that the boundary of $U_a^-$ touches that of $U'_{a,n}$. Thus for $-\pi \leq \psi < 0$ we have $U'_{a,n} \subset U_a^-$. Thus we have shown that for $c \in \RR$, $0<c<t$, $|a|\leq 4$, $n\geq 3$, and $-\pi < \psi \leq \pi$, if $n$ is odd then $U'_{a,n} \subset U_a^-$.
\end{proof}

We apply Douady and Hubbard's techniques to conclude:

\begin{theorem}\label{thm:babyJulia1}
Let $n \geq 3$, $|a| \leq 4$, $\Arg(a) \neq \pi$, $s=\max\{4, |c|, |a| \}$ and $0 < c < \frac{|a|^{1/n}}{s}$. Then for $n$ odd, $\Rnca$ restricted to the set of points whose orbits remain in $U'_{a,n}$
is topologically conjugate to some quadratic polynomial on its filled Julia set,  and $\KRnca$ contains a baby quadratic Julia set.
\end{theorem}

\begin{proof}
We established in Proposition \ref{prop:FirstResult} that with these parameter restrictions we have $U'_{a,n} \subset U_a^-$, and from Proposition~\ref{prop:Uellipse} we know $\Rnca$ maps $U'_{a,n}$ two-to-one onto $U_a^-$.  We need the restriction that $\psi \neq \pi$ so that we can also claim that $U'_{a,n}$ is also relatively compact in $U_a^-$. Then by Definition \ref{defn:DH1}, $\Rnca$ is a polynomial-like map of degree 2 on $U'_{a,n}$ given these parameters. Thus, by Proposition \ref{prop:DH2}, the filled Julia set of $\Rnca$ on $U'_{a,n}$ (the set of points whose orbits remain in $U'_{a,n}$) is topologically conjugate to the filled Julia set of some quadratic polynomial.
\end{proof}

For an example of an $\Rnca$ satisfying this theorem, see Figure \ref{fig:JuliaOverlay1} which shows for the map with $n=3, a=0.25, c=0.2$ both the filled Julia set and an overlay of  $U'_{a,3}$.

%%%%%%%%%%%%%%%%%%%%%%%%%%%%%%%%%%%
\subsection{Results for $|c| > 1 + \epsilon$}
\label{subsection:initialcresults}
%--------------------------------------

Now we turn to the $|c| > 1$ case. First, we provide a polar rectangle containing the boundedness locus $\Mn(\Rnca)$ in an $a$-plane for fixed $c$ with $|c| > 1 + \epsilon$, where $\epsilon > 0$, and $n$ sufficiently large. Using this region, we show $\Rnca$ has a subset on which it is polynomial-like, hence contains a baby quadratic Julia set.
We do this by first establishing for real values of $c$, and then extending the result to complex values of $c$.

\begin{prop} \label{prop:MAATalkProp1}
Given any $\epsilon>0$, $c \in \mathbb{R}$, and $c > 1+\epsilon$ there exists an $N \in \mathbb{N}$ such that for $n \geq N$, $\Mn(\Rnca)$ in the $a$-plane is contained in the polar rectangle $L_c$, where $L_c$ is given by the set of $a$-values such that
$$\frac{(c-(1+\epsilon))^2}{4} \leq |a| \leq \frac{(c+(1+\epsilon))^2}{4}$$
\begin{center}
and
\end{center}
$$-2\sin^{-1}\left(\frac{1 + \epsilon}{c}\right) \leq \Arg(a) \leq 2\sin^{-1} \left(\frac{1+ \epsilon}{c} \right).$$
\end{prop}
\begin{proof}
From Corollary~\ref{cor:BS1} we have that given any $c \in \mathbb{C}, a \in \mathbb{C}^{*}$, and $\epsilon > 0$ there is an $N \geq 2$ such that for all $n \geq N$, $K(\Rnca) \subset \mathbb{D}(0,1+\epsilon).$  Our assumptions satisfy these conditions, so in order for $a$ to be in the boundedness locus for $\Rnca$, we must have at least one of $| c \pm 2\sqrt{a}| < 1+\epsilon$.

Since $Re(2\sqrt{a}) > 0$, $c \in \mathbb{R}$ and $c > 1+\epsilon$ we have that $Re(c+2\sqrt{a}) > 1+\epsilon$, and so the orbit of $c+2\sqrt{a}$ cannot possibly remain bounded. Thus we can only look at $c-2\sqrt{a}$.
Now, if 
$v_- = c-2\sqrt{a} \in \DD(0, 1+\ep),$
%$(c-2\sqrt{a}) \in \{z: |z|<1 + \epsilon \}$,
 then we must have 
% $-2\sqrt{a} \in   \{z: |z+c|<1+\epsilon\}$.
% Hence 
$\sqrt{a} \in 
\DD(c/2, (1+\epsilon)/2).$ 
% \left\lbrace z: |z-{c}/{2}|<{(1+\epsilon)}/{2}\right\rbrace$.
By squaring this region we can use the result to create a set of bounds to contain $\Mn(\Rnca)$. Since $c$ is real and $c > 1+ \epsilon$, this disk is centered on the real axis and $Re(z) > 0$ for all $z$ inside it. Thus for all of 
$\DD(c/2, (1+\epsilon)/2),$
%$\left\lbrace z: |z-{c}/{2}|<{(1+\epsilon)}/{2}\right\rbrace$, 
we get ${(c-(1+\epsilon))}/{2} \leq |\sqrt{a}| \leq ({c+(1+\epsilon)})/{2}$. Squaring this inequality yields the bounds for $|a|$ given in the proposition.
%:${(c-(1+\epsilon))^2}/{4} \leq |a| \leq {(c+(1+\epsilon))^2}/{4}.$

We find bounds for $\Arg(a)$ by finding the largest and smallest arguments that occur on the disk of radius ${(1+\epsilon)}/{2}$ centered at ${c}/{2}$, and then doubling them. Since the disk is centered on the real axis, the largest and smallest arguments will simply be opposites of each other. The largest argument occurs when the line segment extending from the origin is tangent to our disk. We use the right triangle this creates to find the angle. See Figure \ref{fig:Triangle}.

\begin{figure}
\centering
\includegraphics[scale=0.5]{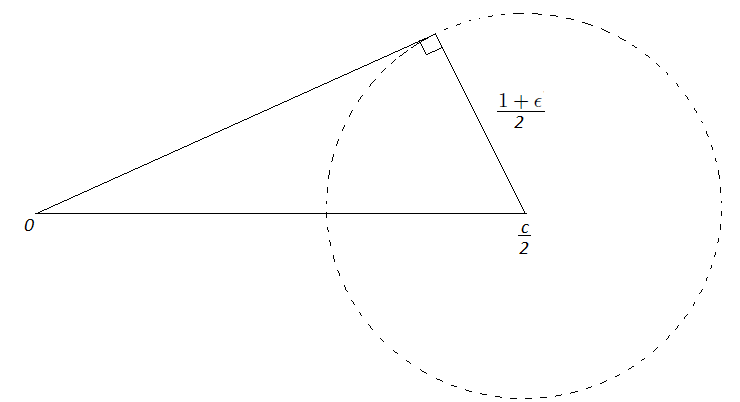}
\caption{Illustration of the largest possible argument of $a$}
\label{fig:Triangle}
\end{figure}

From this we see that the sine of this angle is ${(1+\epsilon)}/{c}$, so the largest argument is $\sin^{-1}({(1+\epsilon)}/{c})$, and the smallest is $-\sin^{-1}({(1+\epsilon)}/{c})$. Doubling these gives us the upper and lower bounds for the argument of $a$ stated in the proposition.
%$$-2\sin^{-1}\left({(1+\epsilon)}/{c}\right) \leq \Arg(a) \leq 2\sin^{-1} \left({(1+\epsilon)}/{c} \right).$$

So, for $n$ sufficiently large, the modulus and argument of $a$ must be bounded in both these ways in order to have $c-2\sqrt{a} \in K(\Rnca)$. 
\end{proof}

We extend the result to complex $c$ via a rotation argument.

\begin{prop} \label{prop:MAATalkProp2}
Given any $\epsilon>0$, $c \in \mathbb{C}$, and $|c| > 1+\epsilon$ there exists an $N \in \mathbb{N}$ such that for $n \geq N$ the slice of $\Mn(\Rnca)$ in the $a$-plane is contained in the polar rectangle $L_c$, where $L_c$ is given by the set of $a$-values such that
$$\frac{(|c|-(1+\epsilon))^2}{4} \leq |a| \leq \frac{(|c|+(1+\epsilon))^2}{4}$$
\begin{center}
and
\end{center}
$$-2\sin^{-1}\left(\frac{1 + \epsilon}{|c|}\right) + 2\Arg(c) \leq \Arg(a) \leq 2\sin^{-1} \left(\frac{1+ \epsilon}{|c|} \right) + 2\Arg(c).$$
\end{prop}
\begin{proof}
The polar rectangle $L_c$ here is simply the polar rectangle from Proposition \ref{prop:MAATalkProp1}, rotated by twice the argument of $c$. Suppose $\left| |c|-2\sqrt{a}\right| < 1+ \epsilon$. Then any rotation of $(|c|-2\sqrt{a})$ will still satisfy the same inequality. Thus,
\begin{align*}
& \left| |c|-2\sqrt{a} \right|  < 1 + \epsilon 
\ \Rightarrow \ 
\left| |c|-2\sqrt{a} \right|  \left| e^{i\Arg(c)} \right|  < 1+ \epsilon 
 \\
& \ \Rightarrow \ 
 \left||c|e^{i\Arg(c)}-2e^{i\Arg(c)}\sqrt{a} \right|  < 1+ \epsilon 
\ \Rightarrow  \left| c-2\sqrt{ae^{2i\Arg(c)}}\right| < 1+ \epsilon. 
\end{align*}
So, if $a$ is rotated by twice $\Arg(c)$, we get $|c-2\sqrt{a}| < 1 + \epsilon$, and therefore this $a$ is potentially in the boundedness locus. We can do this with any $a$ from the polar rectangle from the previous proposition, so we create the new polar rectangle by rotating the previous one by $2\Arg(c)$.
\end{proof}

Now that we have found a polar rectangle $L_c$ containing the boundedness locus, we use this to help show some $U'_{a,k} \subset U_a$ for the values of $a \in L_c$. We again begin by showing $U'_{a,k} \subset \EE$ for all $k$.

\begin{prop} \label{prop:UprimeInEllipse}
Suppose  $c \in \mathbb{C}$, $|c|>1+ \epsilon$, $\epsilon > 0$, and $a \in L_c$ as given in Proposition \ref{prop:MAATalkProp2}. There is an $N \in \mathbb{N}$ such that for all $n \in \mathbb{N}$, $U'_{a,k} \subset \mathcal{E}$ for all $k\in \{0,\ldots,2n-1\}$.
\end{prop}
\begin{proof}
By Proposition \ref{prop:MAATalkProp2}, let $N_1$ be sufficiently large so that $\Mn(\Rnca)$ in this $a$-plane is contained in $L_c$. For a chosen $c$, let $N_2$ be so large that $6|c|+4 < 2^{{N_2}+1}$. Finally, let $N=\max\{N_1,N_2\}$. Now, for all $n \geq N$, we can use the inequalities $6|c|+4 < 2^{n+1}$ and $|a| \leq {(|c|+(1+ \epsilon))^2}/{4}$.

In order to show that $U'_{a,k} \subset \mathcal{E}$ we first show that $U'_{a,k} \subset \mathbb{D}(0,2)$, and then we show $\mathbb{D}(0,2)\subset \mathcal{E}$. Recall the definition of $U'_{a,k}$:
(Definition~\ref{defn:U_ak_prime}):
for $\psi = \Arg(a)$, 
 $$U'_{a,k} = \left\lbrace z \;\middle|\; \frac{|a|^{\frac{1}{n}}}{2} < |z| < 2, \; 
\left| \Arg(z) - \left( \frac{\psi+2k\pi}{2n} \right) \right| < \frac{\pi}{2n} 
\right\rbrace.$$
By our assumptions we have
\begin{equation*}
|a| 
 \leq \frac{(|c|+(1+ \epsilon))^2}{4} 
< \frac{(|c|+|c|))^2}{4}  
= \frac{4|c|^2}{4} 
= |c|^2. 
\end{equation*}
Since $|c|^2 > |a|$ we can write $|c| > |a|^{1/2}$. From this we get:
$$
2^{n+1}  > 6|c|+4 > 6|c| > 2|c| > 2|a|^{1/2}.
$$

Manipulating $2|a|^{1/2} < 2^{n+1}$, we get
\begin{equation*}
2|a|^{1/2}  < 2^{n+1} \Rightarrow \ 
|a|^{1/2}  < 2^n \Rightarrow \ 
|a|  < 4^n \Rightarrow  \ 
|a|^{1/n}  < 4 \Rightarrow \ 
\frac{|a|^{1/n}}{2}  < 2. 
\end{equation*}
This verifies that our definition of $U'_{a,k}$ works nicely for these $a$-values and shows that $U'_{a,k} \subset \mathbb{D}(0,2)$ for every $k$. Next, to show $\mathbb{D}(0,2) \subset \mathcal{E}$, we show that $|z-c+2\sqrt{a}| + |z-c-2\sqrt{a}| < 2^{n+1}+\frac{|a|}{2^{n-1}}$ for all $z \in \mathbb{D}(0,2)$. Let $|z|<2$. Then,
\begin{align*}
&|z-c+2\sqrt{a}| + |z-c-2\sqrt{a}|  \leq 2|z| + 2|c| + 4|a|^{1/2} 
 < 4 + 2|c| + 4|a|^{1/2} \\
& < 4 + 2|c| + 4|c| 
 = 6|c| + 4 
 < 2^{n+1} 
 < 2^{n+1} + \frac{|a|}{2^{n-1}}.
\end{align*}
Thus for sufficiently large $n$ we have $U'_{a,k} \subset \EE$ for all $k$.
\end{proof}

Proposition~\ref{prop:UprimeInEllipse} gives us that $U'_{a,k} \subset \EE$ for all $k\in \{0,\ldots,2n-1\}$ and
 $c \in \mathbb{C}$ with $|c|>1+ \epsilon$, $\epsilon > 0$, $a \in L_c$ as given in Proposition \ref{prop:MAATalkProp2}, and $n \in \mathbb{N}$ sufficiently large.

 Now we want to further restrict the parameters to get $U'_{a,k}$ for some $k$ contained in either $U_a^{\pm}$, but want want to point out there are multiple choices to try. In the next subsection (Section~\ref{subsec:cnegative}), we consider the case of $c$ real and negative. Then in the subsequent subsection (Section~\ref{subsection:largec}), we will instead make a different choice and examine $c$ complex and larger. 
 
%-----------------------------------------------------
\subsection{Baby Julia Sets in the case $c<0$}
\label{subsec:cnegative}
%-----------------------------------------------------

Now, somewhat similar to Proposition~\ref{prop:FirstResult} but now for $c< -1$ we will show $U'_{a,0} \subset U_a^+$ by showing that 
$U'_{a,0}$ intersects the correct half of the ellipse, and does not intersect the minor axis of $\EE$. Since $\EE$ is centered at $c$, we find when we restrict to $c$ real and negative, we can ensure that the minor axis never touches $U'_{a,0}$. 

\begin{prop} \label{prop:UprimeInUforNnegative}
Suppose $c \in \mathbb{R}$, $c < -1 - \epsilon$, and $n> N$ where $N$ is as given by Proposition \ref{prop:UprimeInEllipse}. Then for any $\Rnca$ with $n$, $c$, and $a$ satisfying these conditions we have $U'_{a,0} \subset U_a^+$.
\end{prop}

\begin{proof}
By Proposition \ref{prop:UprimeInEllipse} we have that $U'_{a,0} \subset \EE$. We know that the principle root $w_0=a^{1/2n}$ maps to $v_+=c+2\sqrt{a}$ (see Lemma~\ref{lem:critptsmapping}).

Recall $U'_{a,0}$ is:
$$U'_{a,0} = \left\lbrace z \;\middle|\; \frac{|a|^{\frac{1}{n}}}{2} < |z| < 2, \; \frac{\psi-\pi}{2n} < \Arg(z) < \frac{\psi+\pi}{2n} \right\rbrace.$$
Since $-\pi \leq \psi \leq \pi$ and $n \geq 3$, we have that $\frac{\psi + \pi}{2n} \leq \frac{2\pi}{2n} \leq \frac{\pi}{3}$ and $\frac{\psi - \pi}{2n} \geq -\frac{2\pi}{2n} \geq -\frac{\pi}{3}$. Thus for all $z \in U'_{a,0}$ we have $-\frac{\pi}{3} \leq \Arg(z) \leq \frac{\pi}{3}$. 

Thus for any $a$, since $n\geq 3,$ we have that $U_{a,0}'$ lies in the polar rectangle 
$\{ z \  | \  \frac{|a|^{\frac{1}{n}}}{2} < |z| < 2, \; | \Arg(z) |< \frac{\pi}{3} \}$.

Since $c$ is real, if we select a positive real value for $a$, then $\EE$ is centered on the real axis with major axis along $\RR$. The center of $\EE$ is $c<0$ so the minor axis is parallel to the imaginary axis and lies in the left half plane of $\CC$. Thus $U_a^-$ lies entirely in the left half plane, so is disjoint from $U_{a,0}'$, and since $U_{a,0}' \subset \EE,$ we have $U_{a,0}'\subset U_a^+$ for $a>0$.

Now we consider varying $a$. Since we know $U_{a,0}'\subset U_a^+$ for $a>0$, if we prove that the minor axis of $\EE$ never intersects $U'_{a,0}$ for any $n$, $c$, and $a$ with our restrictions, then we will have $U'_{a,0} \subset U_a^+$ for all our values of $a$.

Since $c$ is real and negative, as $\mathcal{E}$ rotates based on $a$ the first points where the minor axis potentially touches the polar rectangle $\{ z \  | \  \frac{|a|^{\frac{1}{n}}}{2} < |z| < 2, \; | \Arg(z) |< \frac{\pi}{3} \}$
will be at $2e^{\frac{\pi}{3}i}$ and $2e^{-\frac{\pi}{3}i}$,  or $1+\sqrt{3}i$ and $1-\sqrt{3}i$. To find the angle of rotation at which this intersection occurs, we construct the triangle in Figure \ref{fig:UprimeinUTriangle}.

\begin{figure}
\centering
\includegraphics[scale=0.475]{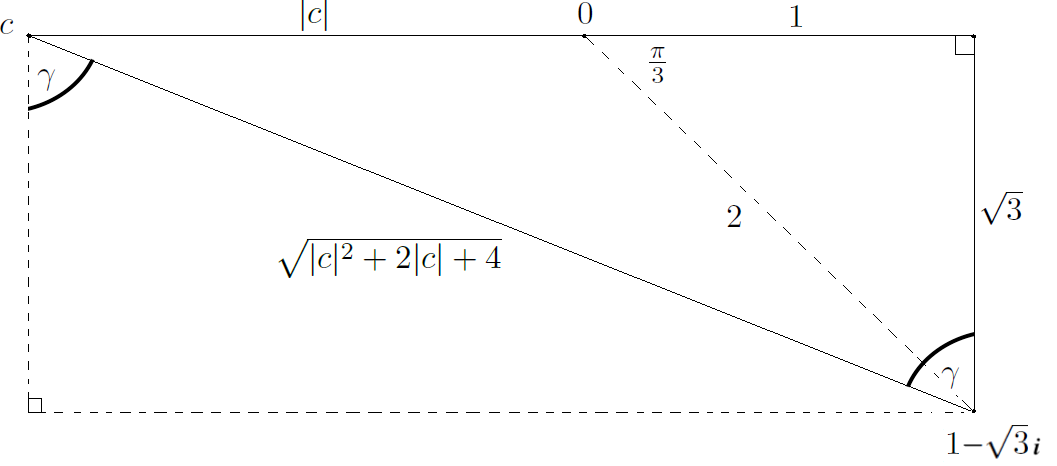}
\caption{Finding the angle of the intersection}
\label{fig:UprimeinUTriangle}
\end{figure}

From this triangle we can determine that $\gamma = \sin^{-1}\left(\frac{|c|+1}{\sqrt{|c|^2+2|c|+4}}\right)$. Thus the minor axis of $\mathcal{E}$ intersects this polar region when the angle of rotation is this $\gamma$.
The angle of rotation of $\mathcal{E}$ is given by ${\psi}/{2}$, so we shall show that given our restrictions we end up with ${\psi}/{2} < \gamma$ for all allowed parameters.

Our restrictions on $\Arg(a)$ are
$$2\Arg(c) - \sin^{-1}\left(\frac{1+ \epsilon}{|c|}\right) \leq \psi \leq 2\Arg(c) + \sin^{-1}\left(\frac{1+ \epsilon}{|c|}\right).$$
Since $\Arg(c)=\pi$, this means that $\psi$ is within $\pm \sin^{-1}\left(\frac{1+ \epsilon}{|c|}\right)$ of $2\pi$, or $0$. So we can say
$$-\sin^{-1}\left(\frac{1+ \epsilon}{|c|}\right) \leq \psi \leq \sin^{-1}\left(\frac{1+ \epsilon}{|c|}\right).$$
Thus by our restrictions we have
$$\frac{\psi}{2} \leq \frac{1}{2}\sin^{-1}\left(\frac{1+ \epsilon}{|c|}\right) < \frac{1}{2}\sin^{-1}(1) = \frac{\pi}{4}.$$

Now we take a look at the contents of our other inverse sine.
Observe from Figure~\ref{fig:UprimeinUTriangle} that when $|c|=1$, the top of the rectangle is length $2$, with the side length $\sqrt{3}$ yields a hypotenuse of length $7$. Thus when $|c|=1$, we have $\sin(\gamma) = 2/\sqrt{7}$, so $\gamma$ is about $.86$ radians or $49^{\circ}$. For $|c|>1$, the angle $\gamma$ is larger but remains less than $\pi/2$ radians or $90^{\circ}$, hence $\sin(\gamma)$ strictly increases. Thus since $|c| > 1$, we have $\sin(\gamma) > 2\sqrt{7},$ and hence $\sin(\gamma) > 2/\sqrt{2}.$
%
%\marginpar{Referee comment 20 suggests a different calculation}
% \begin{align*}
% & \sin (\gamma) = 
%  \frac{|c|+1}{\sqrt{|c|^2+2|c|+4}}  = \frac{\sqrt{(|c|+1)^2}}{\sqrt{|c|^2+2|c|+4}} 
%  = \sqrt{\frac{|c|^2+2|c|+1}{|c|^2+2|c|+4}} \\
% & = \sqrt{\frac{|c|^2+2|c|+4-3}{|c|^2+2|c|+4}} 
%  = \sqrt{1 - \frac{3}{|c|^2+2|c|+4}}. 
% \end{align*}
% Since $|c|>1$, we have $|c|^2+2|c|+4 > 7$, and so $1/(|c|^2+2|c|+4) < 1/7$. Thus,
% $\sin(\gamma)  > \sqrt{1 - {3}/{7}} = {2}/{\sqrt{7}} > {\sqrt{2}}/{2}.$ 
%
Thus 
$ \gamma
> \frac{\pi}{4}$. 
We have now shown that
$$
\frac{\psi}{2} \leq \frac{1}{2}\sin^{-1}\left(\frac{1+ \epsilon}{|c|}\right) < \frac{1}{2}\sin^{-1}(1) = \frac{\pi}{4} < 
\gamma.
$$
A similar argument using an almost identical triangle gives us the following inequality:
$$\frac{\psi}{2} \geq -\frac{1}{2}\sin^{-1}\left(\frac{1+ \epsilon}{|c|}\right) > -\frac{1}{2}\sin^{-1}(1) = -\frac{\pi}{4} > 
- \gamma.
$$

From this we see that the minor axis of $\mathcal{E}$ never intersects the polar rectangle given by $\{ z \  | \  \frac{|a|^{\frac{1}{n}}}{2} < |z| < 2, \; | \Arg(z) |< \frac{\pi}{3} \}$
 for any allowable choice of $n$, $c$, and $a$, and therefore never intersects $U'_{a,0}$ which is contained within. Thus we have $U'_{a,0} \subset U_a^+$ for $c \in \mathbb{R}$, $c < - 1 - \epsilon$, and $n$ sufficiently large.
\end{proof}

Just as with the small $c$ case, we can now conclude that the sets $U'_{a,0}$ with the above restrictions also contain baby Julia sets of quadratic polynomials.

\begin{theorem}\label{thm:babyJulia2}
Suppose $c \in \mathbb{R}$, $c < - 1 - \epsilon$, where $\epsilon>0$, and $n$ sufficiently large so as to satisfy Proposition \ref{prop:UprimeInEllipse}. Given these restrictions, $\Rnca$ restricted to the set of points whose orbits remain in $U'_{a,0}$ is topologically conjugate to some quadratic polynomial on its filled Julia set.
\end{theorem}
\begin{proof}
The proof for this theorem is essentially the same as that of Theorem \ref{thm:babyJulia1}. The result follows from Proposition \ref{prop:DH2}.
\end{proof}

%%%%%%%%%%%%%%%%%%%%%%%%%%%%%%%%%%%%%%%%%%%%%%%%%%%%%%%%%
\subsection{Baby Julia and Mandelbrot sets for larger $c$}
\label{subsection:largec}
%-----------------------------------------------

In this subsection, we follow a different path continuing from subsection~\ref{subsection:initialcresults} than was taken in the prior subsection~\ref{subsec:cnegative}. Instead of restricting to $c$ real and negative, we consider $c$ complex and larger $|c|\geq 6$.  In this case we are able to locate baby $\mandel$'s in parameter planes for $\Rnca$, for $n$ sufficiently large,
and in this section we prove Theorem~\ref{thm:Main_Theorem_APlane}.

\begin{figure}
\centering
\begin{subfigure}{.95\textwidth}
  \centering
  \captionsetup{width=.5\textwidth}
\includegraphics[width=.75\textwidth,keepaspectratio]{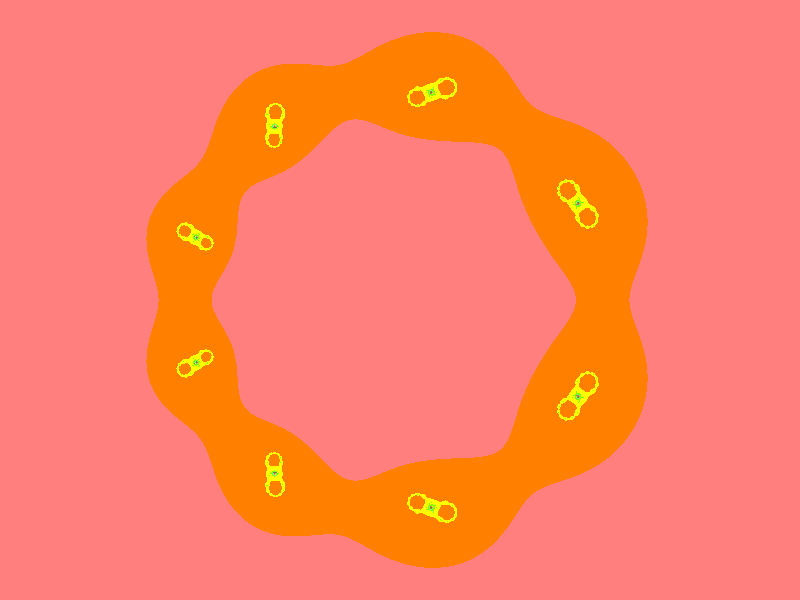}
  \caption{$a$-plane for $n=8$ and $c=6$}
  \label{n8c6}
\end{subfigure}
\begin{subfigure}{.95\textwidth}
  \centering
  \captionsetup{width=.5\textwidth}
\includegraphics[width=.75\textwidth,keepaspectratio]{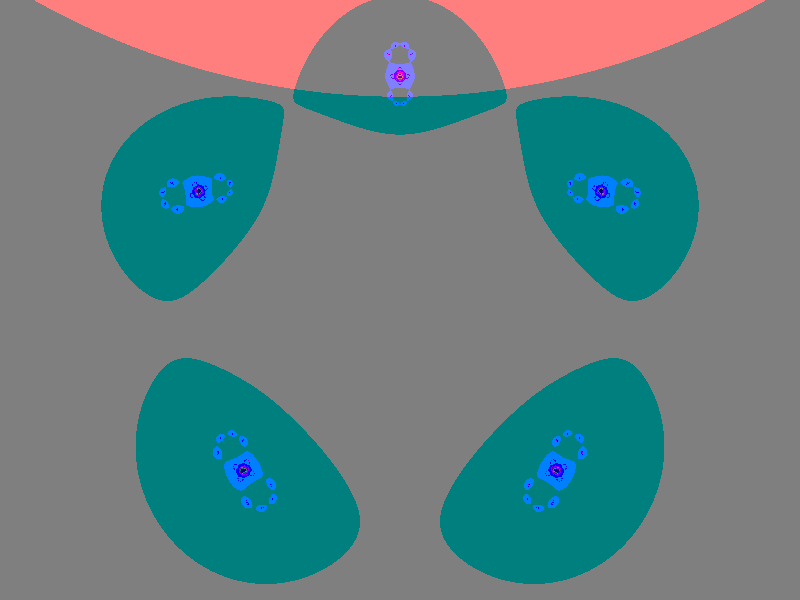}
  \caption{$a$-plane for $n=5$ and $c=-6+6i$}
  \label{n5c-6+6i}
\end{subfigure}
\caption{Two $a$-parameter planes of $\Rnca$ for $|c|\geq 6.$}
\label{fig:RncaParamPlanes}
\end{figure}

Each of the two images in Figure \ref{fig:RncaParamPlanes} will reveal $n$ baby $\mandel$'s in each cluster upon magnification and refinement of the image. 

The purpose of finding the ``spine'' of a boundedness locus is to find a curve within the boundedness locus on which the centers of our primary baby $\mandel$'s lie. This, can allow us to locate the potential baby $\mandel$'s so that we can create appropriate regions to contain them. However, in this section we won't need that method, as instead we'll locate potential centers of the baby $\mandel$'s and use them to define the regions in parameter space which will contain the baby $\mandel$s.

To find these centers, first note by Lemma~\ref{cor:BS2} we have that given any $|c|>1$, for $n \geq 2$, there are $n$ solutions 
to the equation $c \pm 2\sqrt{a} = a^{1/2n}$, i.e. $\Rnca$ maps a critical point to itself. We need some notation to keep track of the dynamics at/near these values.

\begin{defn}
\label{defn:Wparams}
    Let $n\geq 2$ and $c\in \CC$ with $|c|>1$. 
    \begin{itemize}
        \item Set $a_1,\ldots,a_n$ to be the $n$ values of $a$ such that the map  $R_{n,a_j,c}$ has a fixed critical point,
        \item let $w_j$ be the fixed critical point of $a_j$,
        \item let $\psi_j=\Arg(a_j)$,
        \item let $k=k(j)$ be the unique integer in $\{0,\ldots,2n-1\}$ such that \\ $w_j=\xi_k = |a_j|^{1/2n}e^{i \frac{\psi_j +2k\pi}{2n}},$
        \item and for $a$ near $a_j$ (the specific region about $a_j$ is given in Definition~\ref{defn:WforlargeC}), we will use  $U_{a,k} := R(U'_{a,k})$, so
        $$
U_{a,k} =  U^+_{a} \ \text{ for $k$ even, or }  U_{a,k} = U^-_{a} \ \text{ for $k$ odd.} \ \ 
$$
    \end{itemize}
\end{defn}

These fixed critical points $w_j$ are found by making the substitution $w^{2n}=a$ and solving $w^n=({w}-{c})/{2}$. So $w_j^{2n}=a_j$ for each $j=1,\ldots,n$. See Lemma~\ref{lem:critptsmapping} where $\xi_k$ is defined. 
Lemma~\ref{cor:BS2} follows from the proof of the last proposition of \cite{BoydSchulz}. Examining that proof, one can see that since we are in the $|c|>1$ case, there is an annulus containing the critical points, and the fixed critical point $w_j$ with argument $(\psi_j +2k\pi)/2n$ must be in the sector with ray boundaries $\Theta_{n,j} =(\Arg(-c) + (2j-1)\pi)/n$ and  $\Theta_{n,j+1} =(\Arg(-c) + (2j+1)\pi)/n$. As a result,  $k=k(j)$ is the unique integer less than $\pm \frac{1}{2}$ away from $\frac{1}{2\pi} \left(  \Arg(-c) - \frac{\psi_j}{2}  \right).$

Now the $a_j$'s should be the centers of the principal baby $\mandel$'s.
For $|c| \geq 6$ and sufficiently large $n$ we will show that for each $j \in \{ 1, \ldots, n\}$, and for $k = k(j) \in \{ 0, \ldots, 2n-1\},$ given by Definition~\ref{defn:Wparams} (so $w_j= |a_j|^{1/2n}e^{i \frac{\Arg(a_j) +2k\pi}{2n}}$),
there is a region in the $a$-parameter plane, $\Wcj$, such that for all $a \in \Wcj$ we have $U'_{a,k} \subset U_{a_k}$, and  that $\Rnca: U'_{a,k} \rightarrow U_{a,k}$ is polynomial-like of degree two on each $\Wcj$. We will show that for each $a$ in the boundary of any $\Wcj$ we have that the corresponding critical value  
is contained in $U_{a,k}-U'_{a,k}$ and as $a$ travels once around the boundary of any $\Wcj$, the corresponding critical value wraps once around the appropriate critical point in $U'_{a,k}$. Hence by Theorem~\ref{thm:DH3} we will conclude the existence of a baby $\mandel$ in each $\Wcj$.

To begin we will first define the regions $\Wcj,$ for $j=1,\ldots,n$. Each of these regions will contain the entirety of one of our proposed baby $\mandel$'s. 

\begin{defn} \label{defn:WforlargeC} 
Let $|c| \geq 1, n\geq 2,$ and $a_j$ for each $ j=1,\ldots,n$ the $a$-parameter values where $\Rnca$ has a fixed critical point.
Define $\Wcj$ for each $j \in \{1,...,n\}$ as the region in the $a$-parameter plane bounded by the curves:

\begin{align*}
\left( \frac{e^{i\theta}}{4} - \frac{c}{2} \right)^2 
\ \ &\text{for} \ \ 0 \leq \theta \leq  2\pi 
\\
\left( {e^{i\theta}} - \frac{c}{2} \right)^2 
\ \  &\text{for} \ \  0 \leq \theta \leq  2\pi 
\\
\left( \frac{r}{2}e^{i\left(\Arg(w_j)+\frac{\pi}{2n}\right)} - \frac{c}{2} \right)^2 
\ \  &\text{for} \ \  \frac{1}{2} \leq r \leq  2
\\
\left( \frac{r}{2}e^{i\left(\Arg(w_j)-\frac{\pi}{2n}\right)} - \frac{c}{2} \right)^2 
\ \  &\text{for} \ \  \frac{1}{2} \leq r \leq  2.
\end{align*}
\end{defn}

The first two curves above, the inner and outer boundaries, are found, respectively, by solving $|c\pm 2\sqrt{a}|=\frac{1}{2}$ and $|c \pm 2\sqrt{a}|=2$ for $a$. These two curves are both close to circular. 
The third and fourth curves, the ``ray''-like boundaries, are found by solving 
$\Arg(c\pm 2\sqrt{a})=\Arg(w_j)+\frac{\pi}{2n}$ and $\Arg(c \pm 2\sqrt{a})=\Arg(w_j)-\frac{\pi}{2n}$ for $a$.
%\marginpar{Referee says phrase in terms of $\Arg(c\pm 2\sqrt{a}) =$ something so we don't have to introduce $x$. }
%$c\pm 2\sqrt{a}=xe^{i\left(\Arg(w_j)+\frac{\pi}{2n}\right)}$ and $c \pm 2\sqrt{a}=xe^{i\left(\Arg(w_j)-\frac{\pi}{2n}\right)}$ for $a$. 
These latter two curves are both arcs of a parabola.

\begin{figure}
\centering
\begin{subfigure}{.95\textwidth}
\includegraphics[width=.5\textwidth,keepaspectratio]{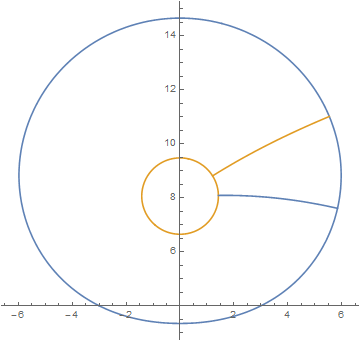} 
\includegraphics[width=.5\textwidth,keepaspectratio]{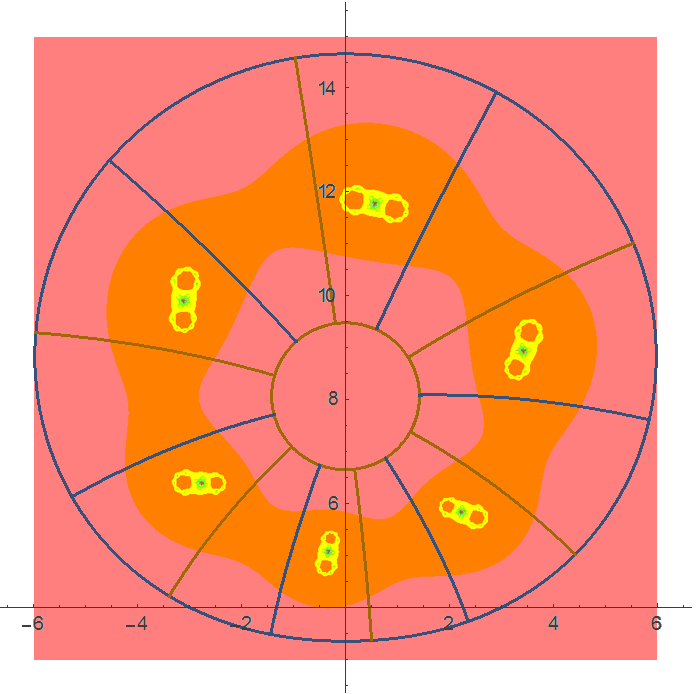} 
\label{fig:SampleW1}
\caption{For $c=4+4i, n=6$, left is a $\Wcj$ and right is an overlay of all $\Wcj$'s}
\end{subfigure}
\centering
\begin{subfigure}{.98\textwidth}
 \includegraphics[width=.5\textwidth,keepaspectratio]{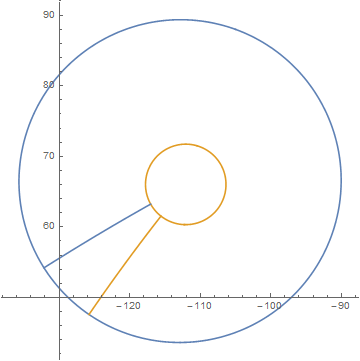} 
\includegraphics[width=.5\textwidth,keepaspectratio]{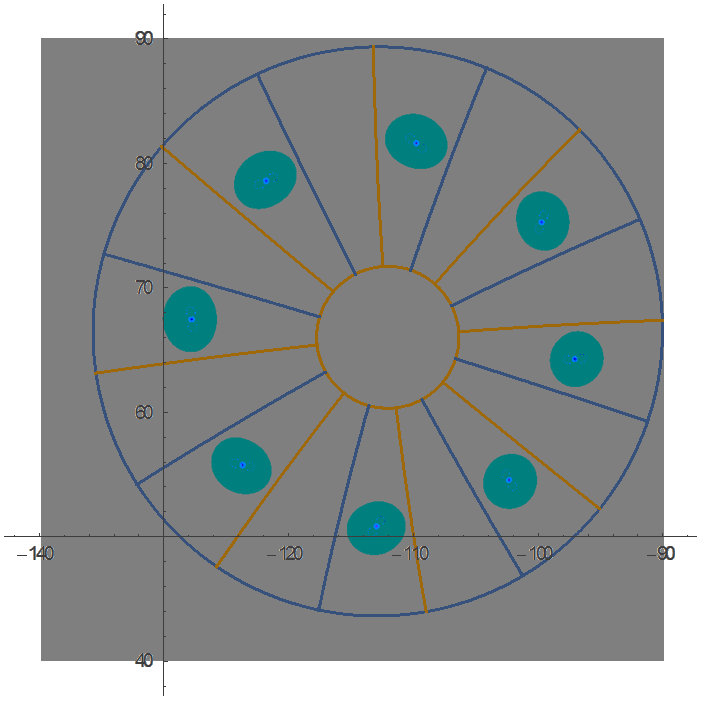} 
\label{fig:SampleW2}
\caption{For  $c=-6-22i, n=8$, left is a $\Wcj$ and right is an overlay of all $\Wcj$'s}
\end{subfigure}
\caption{$\Wcj$'s in the $a$-parameter planes}
\label{fig:SampleW}
\end{figure}

We define $\Wcj$ this way so that (for sufficiently large $c$ and $n$) as the parameter $a$ travels around the boundary of $\Wcj$ the appropriate critical value $R(w_j)$ travels around a specific region, which we will show is in $U_{a,k}-U'_{a,k}$. 
Figure \ref{fig:SampleW} gives an illustration of some $\Wcj$'s.

We want to show that, for all $a$ in each of the $\Wcj$'s, $\Rnca$ is polynomial-like of degree 2 on $U'_{a,k}$.
Recall the definition of $U'_{a,k}$:
$$U'_{a,k} = \left\lbrace z \;\middle|\; \frac{|a|^{\frac{1}{n}}}{2} < |z| < 2, \; 
\left| \Arg(z) -  \left(  \frac{\psi + 2k\pi}{2n} \right)  \right|  < \ \frac{\pi}{2n}.
\right\rbrace
$$

We want to show  that $U'_{a,k(j)} \subset U_{a,k(j)}$ for every $a$ in every $\Wcj$. 
As before, we first establish that each $U'_{a,k}$ is contained within the ellipse $\EE$.

\begin{lem} \label{lem:U'inellipse}
Assume $n$ satisfies $4|c|+8 < 2^{n+1}$.
For any $a$ in any $\Wcj$, we have $\overline{\mathbb{D}(0,2)} \subset \EE$, hence $U'_{a,k} \subset \EE$ (for any $k$).
\end{lem}
\begin{proof}
As specified in Lemma~\ref{prop:Uellipse}, the length of the semi-major axis of $\EE$ is given by $2^n + \frac{|a|}{2^n}$ and the length of the semi-minor axis of $\EE$ is given by $2^n - \frac{|a|}{2^n}$. Recall  $\EE$ is centered at $c$ and is rotated by a degree of ${\psi}/{2}$, and the two foci of $\EE$ occur at $c \pm 2\sqrt{a}$, the two critical values of $\Rnca$. Hence, the distance from the center of the ellipse to either focus is $2\sqrt{|a|}$. 

Now that we know the locations of the two foci of $\EE$ we can prove that $\overline{\mathbb{D}(0,2)} \subset \EE$ by examining the sum $|z - (c+2\sqrt{a})| + |z - (c-2\sqrt{a})|$ for $|z| \leq 2$. The sum of the distances from any point on $\EE$ to each foci is $2 \left( 2^n + \frac{|a|}{2^n} \right)$. If we can show  $|z - (c+2\sqrt{a})| + |z - (c-2\sqrt{a})| < 2 \left( 2^n + \frac{|a|}{2^n} \right)$ for all allowable choices of $a$, this would give us that $\overline{\mathbb{D}(0,2)} \subset \EE$. To do this we first find an upper bound for the magnitude of $a$. 
Recall that for each $\Wcj$ the outer boundary is given by
$$
\left( {e^{i\theta}} - \frac{c}{2} \right)^2 =
e^{2i\theta} - ce^{i\theta} + \frac{c^2}{4}  \;\; \text{for} \;\; 0 \leq \theta \leq 2\pi.
$$
From this we get
\begin{equation*}
|a| \leq \left| e^{2i\theta} - ce^{i\theta} + \frac{c^2}{4} \right|
\leq |e^{2i\theta}| + |ce^{i\theta}| + \left|\frac{c^2}{4}\right|
= \frac{|c|^2}{4} + |c| + 1
= \left( \frac{|c|}{2} + 1 \right)^2.
\end{equation*}
Going back to $|z - (c+2\sqrt{a})| + |z - (c-2\sqrt{a})|$ and using our initial assumption that $n$ is sufficiently large so that $4|c|+8 \leq 2^{n+1}$ we now have
\begin{align*}
& |z - (c+2\sqrt{a})| + |z - (c-2\sqrt{a})| 
\leq 2|z| + 2|c| + 4|a|^{1/2}
\\ &
\leq 4 + 2|c| + 4\left(\frac{|c|}{2} + 1 \right)
= 4|c| + 8
\leq 2^{n+1}
< 2^{n+1} + \frac{|a|}{2^{n-1}}
= 2 \left( 2^n + \frac{|a|}{2^n} \right).
\end{align*}
Thus $\overline{\mathbb{D}(0,2)} \subset \EE$ and therefore $U'_{a,k} \subset \EE$ (for any $k$).
\end{proof}

Now, to get to $U'_{a,k}$ is contained in its image for every $a$ in $W_{c,j}$, we first examine the change in $\Arg(a)$ in the region $W_{c,j}$.

\begin{lem} \label{lem:Critvaluearg}
Assume $|c| \geq 6$ and $n\geq 3$.
As $a$ varies throughout any given $\Wcj$, the argument of any critical point $a^{1/2n}$ changes by less than $\frac{\pi}{2n}$ .
\end{lem}
\begin{proof}
We mentioned earlier that as we trace $a$ around the boundary of a given $\Wcj$ the associated critical value travels around a specific region in the dynamical plane. We will call this region $\Vcj$:
\begin{equation} \label{eqn:Vcj}
\Vcj= \left\lbrace z \;\middle|\; \frac{1}{2} \leq |z| \leq 2, \;  | \Arg(z) - \Arg(w_j) | \leq \frac{\pi}{2n} \right\rbrace.
\end{equation}
The two ``ray-like'' boundaries of $\Vcj$ are created when $a$ is on the two boundary curves of $\Wcj$ given by 
$$
\left( \frac{r}{2}e^{i\left(\Arg(w_j)\pm \frac{\pi}{2n}\right)} - \frac{c}{2} \right)^2 
$$
for $\frac{1}{2} \leq r \leq 2.$

First we show for any fixed $r$ in $[\frac{1}{2}, 2]$, the change in argument, as $x: -\frac{\pi}{2n} \to \frac{\pi}{2n}$, of
$\left( \frac{r}{2}e^{i\left(\Arg(w_j)+x\right)} - \frac{c}{2} \right)^2 $
 is less than $\frac{2\pi}{n}$.

If we use the fact that $\Arg(re^{i(\theta + x)} + se^{i\psi}) < \Arg(re^{i\theta} + se^{i\psi}) + x$ when $r \neq 0$ and $s \neq 0$, we see that as $x: -\frac{\pi}{2n} \to \frac{\pi}{2n}$, $\Arg\left( \frac{r}{2}e^{i(\Arg(w_j)+x)} - \frac{c}{2} \right)$ changes by less than $2(\frac{\pi}{2n}) = \frac{\pi}{n}$.
Thus, as $x: -\frac{\pi}{2n} \to \frac{\pi}{2n}$, $\Arg\left( \frac{r}{2}e^{i(\Arg(w_j)+x)} - \frac{c}{2} \right)^2$ changes by less than $\frac{2\pi}{n}$.
So, in crossing over $\Wcj$ from one parabolic arc to another gives a change in argument of less than $\frac{2\pi}{n}$. If $\Arg(a)$ changes by less than $\frac{2\pi}{n}$, then $\Arg(a^{1/2n})$ for any solution of $a^{1/2n}$ changes by less than $\frac{\pi}{n^2}$. Since $n \geq 3$
we can say $\frac{\pi}{n^2} \leq \frac{\pi}{3n}$.

Next, the two inner and outer boundaries of $\Vcj$ are created when $a$ is on the two boundary curves of $\Wcj$, given by
$$
\left(\frac{e^{i\theta}}{4} - \frac{c}{2}\right)^2
\ \;\; \text{  and  }  \;\;  
 \left(e^{i\theta} - \frac{c}{2}\right)^2
 \;\; \text{for} \;\; 0 \leq \theta \leq 2\pi.
$$
Here we show that, for fixed $\theta$, the change in $\Arg(a)$ as we go from the inner to the outer bound is less than ${\pi}/{6}$. 

Since $c$ and $\theta$ are fixed we can think of $\left( re^{i\theta} - \frac{c}{2} \right)$ as a triangle in $\mathbb{C}$ and examine what happens to the argument as we take $r: \frac{1}{4} \to 1$. 
The change in argument for $\left( re^{i\theta} - \frac{c}{2} \right)$ as $r: \frac{1}{4} \to 1$ will be greatest when $\theta = \Arg(c) \pm \frac{\pi}{2}$. We assume $|c| \geq 6$,
%* \marginpar{* here we used the assumption $|c|>6$} 
so ${|c|}/{2} \geq 3$ and we have a right triangle where the shortest leg is increasing in length from $\frac{1}{4}$ to $1$. Since ${|c|}/{2} \geq 3$ we have $\Arg\left( \frac{e^{i\theta}}{4} - \frac{c}{2} \right) \leq \tan^{-1}({1}/{12})$ and $\Arg\left( e^{i\theta} - \frac{c}{2} \right) \leq \tan^{-1}({1}/{3})$. The change in argument will be greater the smaller $|c|$ is, so the change in argument for $\left( re^{i\theta} - \frac{c}{2} \right)$ as $r: \frac{1}{4} \to 1$ is at most $\tan^{-1}({1}/{3}) - \tan^{-1}({1}/{12}) < {\pi}/{12}$.

The change in argument for $\left( re^{i\theta} - \frac{c}{2} \right)$ as $r: \frac{1}{4} \to 1$ is less than ${\pi}/{12}$, so the change in argument for $\left( re^{i\theta} - \frac{c}{2} \right)^2$ as $r: \frac{1}{4} \to 1$ is less than ${\pi}/{6}$. If $\Arg(a)$ changes by less than ${\pi}/{6}$, then $\Arg({a^{1/2n}})$ changes by less than $\frac{\pi}{12n}$.

As $a$ travels from any point in a given $\Wcj$ to another, we can think of this change as a shift along two arcs, one of each of the two types detailed above. Therefore we can say that the change in argument of $a^{1/2n}$ as $a$ travels from any one point in $\Wcj$ to another is less than $\frac{\pi}{n^2} + \frac{\pi}{12n}$. Now,
$$\frac{\pi}{n^2} + \frac{\pi}{12n} \leq \frac{\pi}{3n} + \frac{\pi}{12n} = \frac{5\pi}{12n} < \frac{6\pi}{12n} = \frac{\pi}{2n},$$
so as $a$ varies through any one $\Wcj$,  $\Arg(a^{1/2n})$ changes by less than $\frac{\pi}{2n}$, for any solution of $a^{1/2n}$.
\end{proof}

Now, we can show $U'_{a,k(j)}$ is contained in its image for every $a$ in $W_{c,j}$.

\begin{lem} \label{lem:U'inU}
Assume $|c| \geq 6$ and $n$ s.t.\ $4|c|+8 < 2^{n+1}$.
For any $a$ in any $\Wcj$, and for $k=k(j)$, as given by Definition~\ref{defn:Wparams}, $U'_{a,k} \subset U_{a,k}$.
\end{lem}

\begin{proof}
By Lemma~\ref{lem:U'inellipse}, since $n$ satisfies $4|c|+8 < 2^{n+1}$, for any $j$ and any $a\in \Wcj$, we have $U'_{a,k} \subset \EE$. 
We know $U'_{a,k}$ varies analytically with $a$ and there exists an $a_j$ in each $\Wcj$ with a fixed critical point, i.e., with a $w_j\in U'_{a_j,k}$ satisfying $w_j^{2n}=a_j$  and $R_{n,a_j,c}(w_j)=w_j$, so $U'_{a_j,k} \cap U_{a_j,k} \neq \emptyset$.
Now we will show that $U'_{a,k}$ never intersects the minor axis of $\EE$. 
To do so we show
$U'_{a,k} \subset \mathbb{D}(c \pm 2\sqrt{a}, 2\sqrt{|a|}),$ for whichever critical value is the appropriate one. This shows that $U'_{a,k}$ avoids the minor axis, since $\EE$ is the ellipse centered at $c$ with foci at $c\pm 2\sqrt{a}.$

Recall in each $\Wcj$, as $a$ travels around the boundary of $\Wcj$ in the parameter plane, the associated critical value travels around the boundary of the region 
%in the dynamical plane 
which we called $\Vcj$ in Equation~\ref{eqn:Vcj}:
$$
\Vcj= \left\lbrace z \;\middle|\; \frac{1}{2} \leq |z| \leq 2, \;  | \Arg(z) - \Arg(w_j) | \leq \frac{\pi}{2n} \right\rbrace,
$$
which follows immediately from the definition of $\Wcj$ (look back at the paragraph following Definition~\ref{defn:WforlargeC} for $a \in \Wcj$). 

Since $\Arg(w_j) = (\psi_j+ 2k\pi)/2n$, where $\psi_j = \Arg(a_j)$, we have
$$
\Vcj= \left\lbrace z \;\middle|\; \frac{1}{2} \leq |z| \leq 2, \;  
\left| \Arg(z) - \left( \frac{\psi_j+ 2k\pi}{2n} \right) \right| \leq \frac{\pi}{2n} \right\rbrace,
$$
Now recall 
$$U'_{a,k} = \left\lbrace z \;\middle|\; \frac{|a|^{\frac{1}{n}}}{2} < |z| < 2, \; 
\left| \Arg(z) -  \left(  \frac{\psi + 2k\pi}{2n} \right)  \right|  < \ \frac{\pi}{2n}
\right\rbrace.
$$
%-----------------
First, note the inner arc boundary of $U'_{a,k}$ has magnitude $\frac{|a|^{1/n}}{2}$ while the inner arc that the critical value traces has magnitude $\frac{1}{2}$. Recall from the definition of $\Wcj$ that the outer boundary of $\Wcj$ is given by $ \left(e^{i\theta} - \frac{c}{2}\right)^2 = e^{2i\theta} + ce^{i\theta} + \frac{c^2}{4}$, hence
\begin{equation*}
|a| \geq \left| \frac{c^2}{4} \right| - |ce^{i\theta}| - |e^{2i\theta}| 
= \frac{|c|^2}{4} - |c| - 1
\geq \frac{(6)^2}{4} - 7 
= 2.
\end{equation*}
Thus $|a| \geq 2 > 1$, and so $\frac{1}{2} < \frac{|a|^{1/n}}{2}$.
%------------------

Now, note that $U'_{a,k}$ is angularly centered around the argument of the $k^{th}$ critical point of $\Rnca$, that is, $(\psi + 2k\pi)/2n,$ and the difference between the arguments of the $k^{th}$ critical points of the map at $a_j$ versus any $a\in \Wcj$ is $| \psi - \psi_j | / 2n,$ where $\psi=\Arg(a).$

We showed in Lemma~\ref{lem:Critvaluearg} that as $a$ varies in $\Wcj$, the argument of any critical point, any $\Arg(a^{1/2n})$, can change by no more than $\frac{\pi}{2n}$, assuming as we have here that $|c|\geq 6$ and $n\geq 3$ (which we do, since note that $4|c|+8 < 2^{n+1}$ implies $n\geq 5$ so we have $n \geq 3$). So $| \psi - \psi_j | / 2n \leq \pi/2n$, and the difference between the angular centers of $\Vcj$ and $U'_{a,k}$ is at most $\pi/2n$.
Since each of $\Vcj$ and $U'_{a,k}$ are of argument width $\pm \pi/2n$ about their centers,  the maximum difference in argument between a point in $\Vcj$ (like, the appropriate critical value) and any point in $U'_{a,k}$ is $3\pi/2n$. 
%Since $n\geq 3,$  we get $\frac{3\pi}{2n} \leq \frac{\pi}{2}$.

Now, the radius of the smallest circle centered at $c \pm 2\sqrt{a}$ (whichever is the correct one, but the worst case scenario is when it is on a corner boundary point of $\Vcj$ away from $U'_{a,k}$) that is guaranteed to contain $U'_{a,k}$ for any $a$ in $\Wcj$ is found using the triangle in Figure \ref{Diskofradius2rootaTriangle}.
\begin{figure}
\centering
\includegraphics[scale=.85]{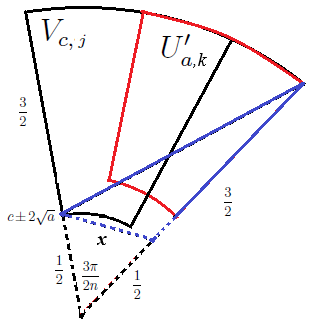}
\caption{Finding the smallest circle which contains $U'_{a,k}$}
\label{Diskofradius2rootaTriangle}
\end{figure}
By the preceding paragraph, the difference in argument between $c \pm 2\sqrt{a}$ and the furthest corner of $U'_{a,k}$ can be no more than $\frac{3\pi}{2n}$, and since $n \geq 3$ we get $\frac{3\pi}{2n} \leq \frac{\pi}{2}$. Thus one of the short legs of the triangle, labeled $x$, has length  $\leq \sqrt{\frac{1}{2}^2 + \frac{1}{2}^2} = \frac{\sqrt{2}}{2}$. The other short leg has length $\frac{3}{2}$. Hence the length of the hypotenuse, equal to  the radius of the smallest circle at $c \pm 2\sqrt{a}$ that is guaranteed to contain $U'_{a,k}$, is $\leq \frac{3+\sqrt{2}}{2}$.
Now, we concluded a couple of paragraphs above that $|a| \geq 2$,
so $2\sqrt{|a|} \geq 2\sqrt{2}$, which is greater than $({3+\sqrt{2}})/{2}$. Thus for any $a \in \Wcj$ a disk of radius $2\sqrt{|a|}$ centered at the appropriate critical value will contain $U'_{a,k}$. This guarantees us that $U'_{a,k} \subset U_{a,k}$.
\end{proof}

And combining the above lemmas, we have shown the following.

\begin{prop} \label{prop:RPolynomialLike}
Let $|c| \geq 6$ and $n$ be large enough that $4|c|+8 < 2^{n+1}$.
The family of functions $\Rnca$ is polynomial-like of degree 2 on each $U'_{a,k(j)}$ for all $a$ in any $\Wcj$. 
\end{prop}

\begin{proof}
%We need $U'_{a,k} \subset U_{a,k}$ and $\Rnca$ maps $U'_{a,k}$ two-to-one onto $U_{a,k}$.
%
The description of $U_{a,k}$ in Lemma~\ref{prop:Uellipse} also applies here as well, so we know that $\Rnca$ maps $U'_{a,k}$ two-to-one onto $U_{a,k}$ (for every $k$). 

So, we just need that $U'_{a,k} \subset U_{a,k}$ for every $a$ in every $\Wcj$, which is established in Lemma~\ref{lem:U'inU}. 

Now we can say that $U'_{a,k} \subset U_{a,k}$ and $U'_{a,k}$ is relatively compact in $U_{a,k}$ for any choice of $a$ in any of the $\Wcj$'s. This combined with the fact that $\Rnca$ is a proper analytic map which maps $U'_{a,k}$ two-to-one onto $U_{a,k}$ gives us that $\Rnca$ is a polynomial-like map of degree two on the orbits which remain in each $U'_{a,k},$ for any $a$ in any $\Wcj$. This ends the proof of the proposition.
\end{proof}

So far we have satisfied the first four conditions in Theorem \ref{thm:DH3} set forth by Douady and Hubbard to show that we have a homeomorphic copy of the Mandelbrot set:
\begin{enumerate}
\item  Each $\Wcj$ is a closed topological disk contained in an open set in $\CC$.
\item The boundaries of $U'_{a,k}$ and $U_{a,k}$ both vary analytically with $a$.
\item The map $(a,z) \mapsto \Rnca(z)$ depends analytically on both $a$ and $z$.
\item On each $\Wcj$, $\Rnca$ restricted to $U_{a,k(j)}'$ has a unique critical point and is polynomial-like of degree two on the orbits which remain in $U_{a,k}'$.
\end{enumerate}

All that remains is to show is that as $a$ winds once around the boundary of any $\Wcj$ the (appropriate) critical value winds once around the (appropriate) critical point, in $U_{a,k} - U'_{a,k}$.
Then we will have that for each $\Wcj$, there is a corresponding baby $\mandel$.

\begin{lem}\label{lem:winding}
Let $|c| \geq 6$ and $n$ be large enough that $4|c|+8 < 2^{n+1}$.
On each region $\Wcj$ the corresponding function given by $\Rnca$ has 
a unique critical point in $U_{a,k}'$.
As we trace the parameter $a$ once around the boundary of any $\Wcj$, the image of this unique critical point winds once around the critical point itself, and remains in $U_{a,k} - U'_{a,k}$.
\end{lem}

\begin{proof}
Recall that each region $\Wcj$ contains a parameter value $a_j$ for which a critical point, $w_j=a_j^{1/2n}$, is fixed, and $\Arg(w_j) = (\psi_j+ 2k\pi)/2n$, where $\psi_j = \Arg(a_j)$.

As mentioned in Lemmas~\ref{lem:Critvaluearg} and~\ref{lem:U'inU}, 
as $a$ winds once around the boundary of any $\Wcj$ the appropriate critical value, $c + 2\sqrt{a}$ or $c-2\sqrt{a}$, winds once around the boundary of
 $\Vcj$:
$$
\Vcj = \left\lbrace z \;\middle|\; \frac{1}{2} \leq |z| \leq 2, \;  
\left| \Arg(z) - \left( \frac{\psi_j+ 2k\pi}{2n} \right) \right| \leq \frac{\pi}{2n} \right\rbrace.
$$
Now recall 
$$U'_{a,k} = \left\lbrace z \;\middle|\; \frac{|a|^{\frac{1}{n}}}{2} < |z| < 2, \; 
\left| \Arg(z) -  \left(  \frac{\psi + 2k\pi}{2n} \right)  \right|  < \ \frac{\pi}{2n} 
\right\rbrace.
$$

We know that each $U'_{a,k}$ contains a critical point, whose image critical value lies in the image of $U'_{a,k}$, called $U_{a,k}$, and that $U'_{a,k} \subset U_{a,k}$, and hence $\partial U_{a,k} \subset U_{a,k} - U'_{a,k}.$ 

In the proof of Lemma~\ref{lem:U'inU} we showed (under the same hypotheses on $c$ and $n$) that $ \frac{1}{2} < \frac{|a|^{\frac{1}{n}}}{2}$. Thus by the definitions above, we know as $a$ traverses the inner boundary of $\Wcj$, so the critical value traverse the inner boundary of $\Vcj$, then the appropriate critical value has modulus $1/2$ and traverses an arc just below the outer arc boundary of $U'_{a,k}$ hence in $U_{a,k} - U'_{a,k}$. As $a$ traverses the outer boundary of $\Wcj$, then the appropriate critical value traverse the outer boundary of $\Vcj$ (of modulus $2$) so traverses the outer boundary of $\partial U_{a,k} \subset U_{a,k} - U'_{a,k}.$

Now we examine the ray boundaries of $\Vcj$.
The two ray boundaries of $\Vcj$ are created when the appropriate critical value is on the two boundary curves of $\Vcj$ so $\Arg(z) = (\psi_j + (2k\pm 1)\pi)/2n$. 
The set $U'_{a,k}$ 
has boundaries of argument $(\psi + (2k\pm 1)\pi)/2n$.

But then, to know that the critical value is in $\partial U'_{a,k} \subset U_{a,k} - U'_{a,k}$ we need only that 
$$
\left|  \frac{\psi_j + (2k \pm 1)\pi}{2n} - \frac{\psi + (2k \pm 1)\pi}{2n}  \right| \geq \frac{\pi}{2n},
$$
which is of course true since $| \psi - \psi_j | \geq 0.$

Hence as $a$ traverses the boundary of $\Wcj$, the critical value wraps around $U_{a,k} - U'_{a,k}$.

\end{proof}

Now we have established our first main theorem:

\begin{proof}[Proof of Theorem~\ref{thm:Main_Theorem_APlane}.]
Proposition \ref{prop:RPolynomialLike} gives us that $\Rnca$ is polynomial-like of degree 2 on each $\Wcj$. This combined with Lemma \ref{lem:winding} finishes our proof. For each $j \in \{1,2,...,n\}$, the Mandelbrot set is homeomorphic to the set of all $a \in \Wcj$ for which the orbit of the appropriate critical point does not escape from $U'_{a,k(j)}$. There are $n$ such parameters $a_j$ and $n$ such regions, so there must exist at least $n$ homeomorphic copies of the Mandelbrot set within the $a$-parameter planes of $\Rnca$.
\end{proof}

%%%%%%%%%%%%%%%%%%%%%%%%%%%%%%%%%%%%%%%%%%%%%%%%
\section{Investigating Diagonal Slices of the $(a,c)$ Parameter Plane}
\label{sec:LinearSlices}

% This section is about 5 1/2 pages. 

\begin{figure}
\centering
\begin{subfigure}{.49\textwidth}
  \centering
 \includegraphics[width=1.0\textwidth,keepaspectratio]{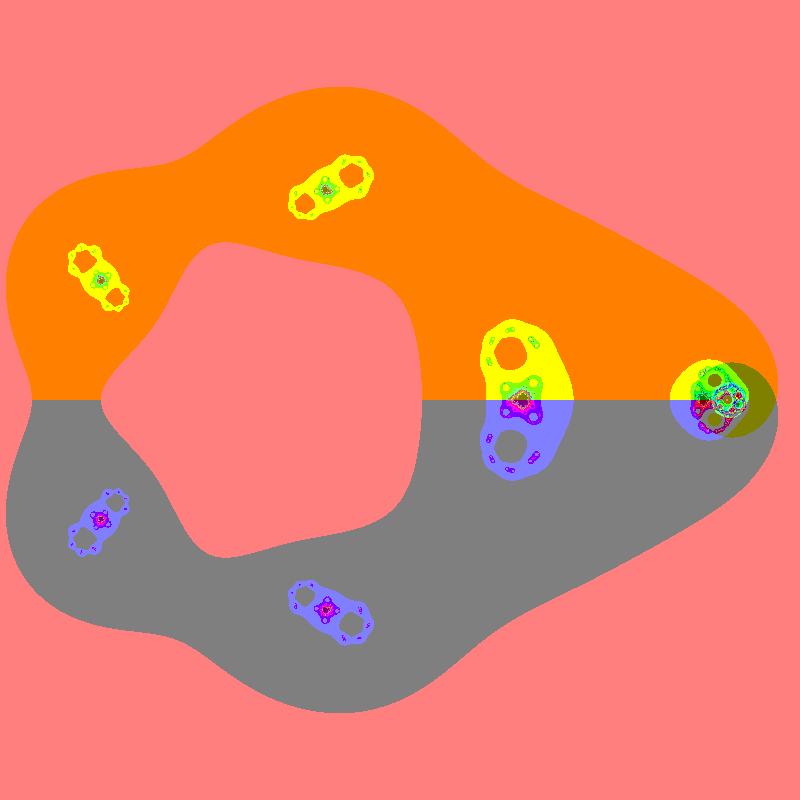} 
  \caption{$t=-0.8i$}
  \label{c=ta_t=-0.8i}
\end{subfigure}
\begin{subfigure}{.49\textwidth}
  \centering
 \includegraphics[width=1.0\textwidth,keepaspectratio]{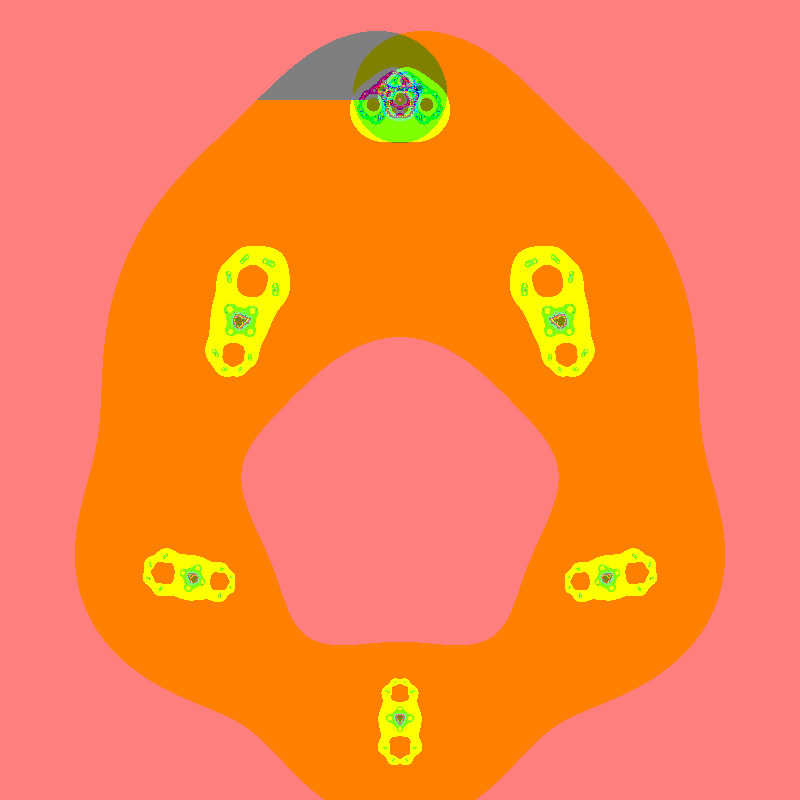} 
  \caption{$t={1}/{\sqrt{2}}+i/{\sqrt{2}}$}
  \label{c=ta_t=1_45degrees}
\end{subfigure}
\caption{Parameter planes of $R_{5,a,ta}$ for some complex $t$ values.}
\label{fig:c=ta}
\end{figure}

%%%%%%%%%%%%%%%%%%%%%%%%%%%%%%%%%%%%%%%%%%%%%%%%
Now we shift to examine  one-dimensional slices with $c=ta,$ 
$$\Rnta(z) = z^n + \frac{a}{z^n} + ta,$$
where $n \geq 3$ fixed, $t \in \CC^*$ fixed, and allow $a \in \CC^*$ to vary analytically. Now $c=ta$ also varies as $a$ varies. Figure \ref{fig:c=ta} shows a few examples of these parameter lines for some complex values of $t$ when $n=5$. Figure~\ref{fig:Spines} shows more examples, along with the spines defined later in this section.

We see again what appear to be small copies of the Mandelbrot set in these parameter planes. Recall Figure \ref{fig:c=a_n=6_intro} shows one such slice and a zoomed-in example of what appears to be a baby Mandelbrot set.

%%%%%%%%%%%%%%%%%%%%%%%%
\subsection{The Spines of $\Rnta$}

For $\Rnta$ we locate the boundedness locus by first finding a spine, in a neighborhood of which we later show the boundedness locus is contained.

 By Corollary~\ref{cor:BS1}, Julia sets for $\Rnca$ are near the unit circle, so that is still the case when $c=ta$. 
So in the subfamily with $c=ta$, filled Julia sets  still fall within the annulus $\mathbb{A}(1-\epsilon, 1+\epsilon)$ for $n$  sufficiently large.

To find the spines, we use that since the Julia sets must be near the unit circle, bounded critical orbits must be near the unit circle as well. 

For any $t \in \CC^*$, let $\mathcal{S}_t$ denote the following set:
\begin{equation}
    \label{eqn:spine}
\mathcal{S}_t = \left\lbrace a=\frac{2}{t^2} + \frac{1}{t}e^{i\theta} \pm \frac{2}{t^2}\sqrt{1+te^{i\theta}} \;\middle|\; 0 \leq \theta \leq 2\pi \right\rbrace.    
\end{equation}

\begin{prop}
\label{prop:critvalspine}
For any $t \in \CC^*$, $\mathcal{S}_t $ is the set of all $a \in \CC$ for which $|ta \pm 2\sqrt{a}| = 1$.
\end{prop}

\begin{proof}
We can write $|ta \pm 2\sqrt{a}| = 1$ as $ta \pm 2\sqrt{a} = e^{i\theta}$ for some $\theta$. To solve this for $a$, we perform a substitution and use the quadratic equation.\\
The solutions to $tu^2 \pm 2u - e^{i\theta} = 0$ are given by
$u = (\mp 1 \pm \sqrt{1 + te^{i\theta}})/{t},$
so the solutions to $ta \pm 2\sqrt{a} - e^{i\theta} = 0$ are 
$$a = \left( \frac{\mp 1 \pm \sqrt{1 + te^{i\theta}}}{t}\right)^2 = \frac{2}{t^2} + \frac{1}{t}e^{i\theta} \pm \frac{2}{t^2}\sqrt{1+te^{i\theta}}.$$
Thus, the collection of $a$-values for which $|ta \pm 2\sqrt{a}| = 1$ is
$$\left\lbrace \frac{2}{t^2} + \frac{1}{t}e^{i\theta} \pm \frac{2}{t^2}\sqrt{1+te^{i\theta}} \;\middle|\; 0 \leq \theta \leq 2\pi \right\rbrace.$$
\end{proof}

\begin{figure}
\centering
\begin{subfigure}{.5\textwidth}
  \centering
  \captionsetup{width=.8\textwidth}
  \includegraphics[width=.95\textwidth,keepaspectratio]{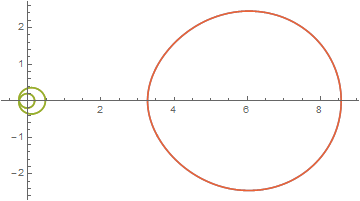}
  \caption{Spine for $t=0.8$}
  \label{Spine_t=0.8}
\end{subfigure}%
\begin{subfigure}{.5\textwidth}
  \centering
  \captionsetup{width=.8\textwidth}
   \includegraphics[width=.95\textwidth,keepaspectratio]{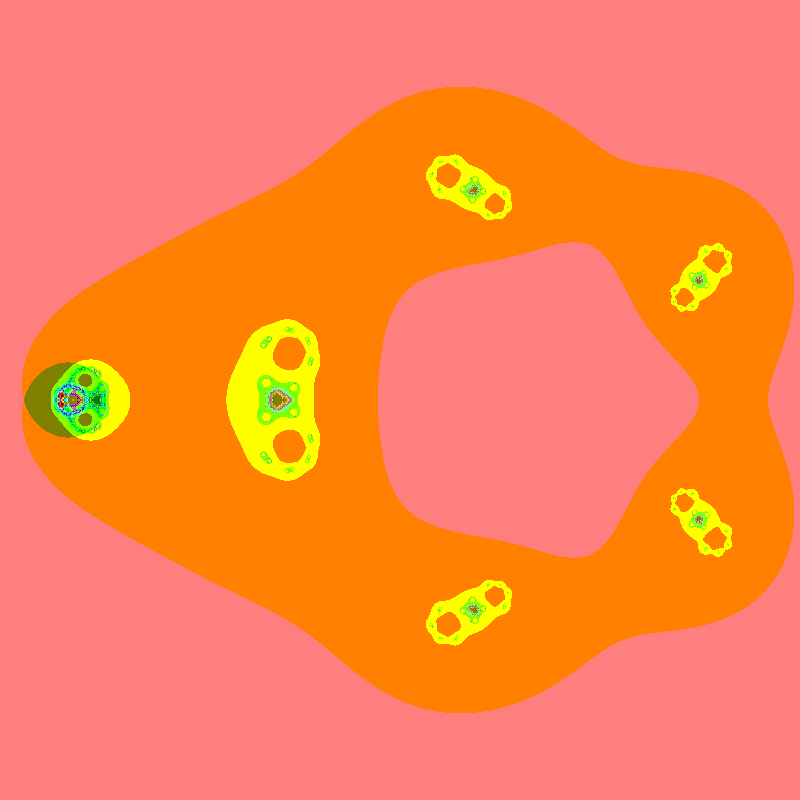} 
  \caption{Parameter plane $t=0.8$}
\end{subfigure}
\begin{subfigure}{.5\textwidth}
  \centering
  \captionsetup{width=.8\textwidth}
  \includegraphics[width=.95\textwidth,keepaspectratio]{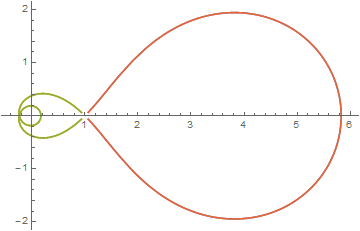}
  \caption{Spine for $t=1$}
  \label{Spine_t=1}
\end{subfigure}%
\begin{subfigure}{.5\textwidth}
  \centering
  \captionsetup{width=.8\textwidth}
  \includegraphics[width=.95\textwidth,keepaspectratio]{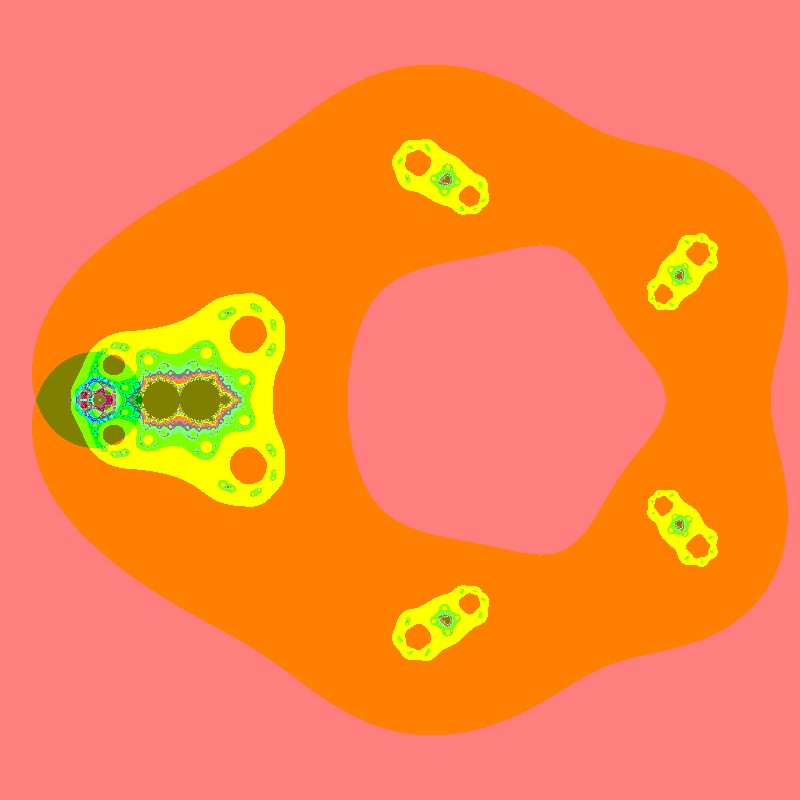} 
  \caption{Parameter plane $t=1$}
\end{subfigure}
\begin{subfigure}{.5\textwidth}
  \centering
  \captionsetup{width=.8\textwidth}
  \includegraphics[width=.95\textwidth,keepaspectratio]{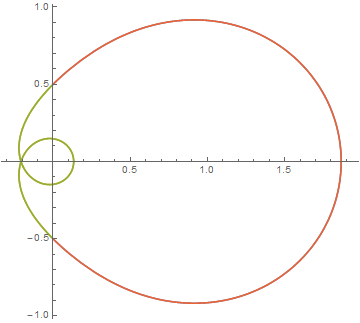}
  \caption{Spine for $t=2$}
  \label{Spine_t=2}
\end{subfigure}%
\begin{subfigure}{.5\textwidth}
  \centering
  \captionsetup{width=.8\textwidth}
   \includegraphics[width=.95\textwidth,keepaspectratio]{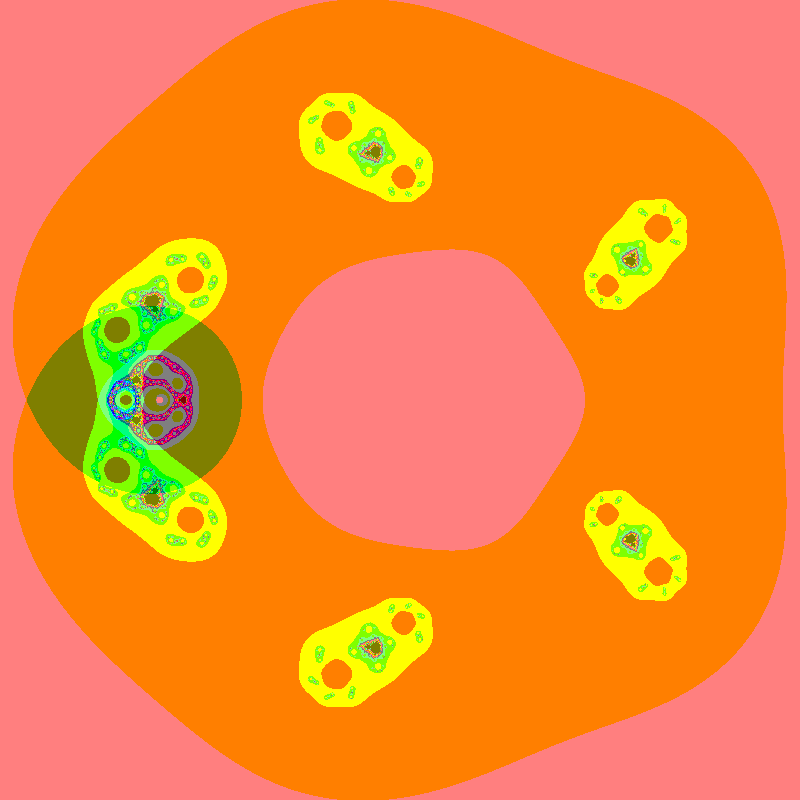} 
  \caption{Parameter plane  $t=2$}
\end{subfigure}
\caption{Spines and the boundedness locus of $R_{5,ta,a}$}
\label{fig:Spines}
\end{figure}

We refer to these curves as the \textit{spines} of the boundedness loci of $\Rnta$.
Figure \ref{fig:Spines} shows three of these curves alongside their corresponding $a$-planes. Note that when $t=1$ we get a bifurcation, where values of $t$ less than one produce a spine which is split into two disconnected curves.

%%%%%%%%%%%%%%%%%%%%%%%%%%%%%%%%%%%%%%%%%%%%
\subsection{The Boundedness Locus of $\Rnta$}

Now we show that $\Mn(\Rnta)$ for each $t\neq 0$ must lie within an $\epsilon$-neighborhood of our spine $\mathcal{S}_t$. We first contain $\Mn(\Rnta)$ in an annulus based on $|t|$. Set
$$u(t) = \frac{2}{|t|^2} + \frac{1}{|t|} + \frac{2}{|t|^2}\sqrt{1+|t|}
\text{ and }
l(t) = \frac{2}{|t|^2} + \frac{1}{|t|} - \frac{2}{|t|^2}\sqrt{1+|t|}.$$

\begin{lem} \label{lem:MandelbrotInAnnulus}
Let $\epsilon > 0$, $t \in \CC^*$. There is an $N \geq 2$ such that for all $n \geq N$, 
$M_n(\Rnta) \subset \mathbb{A}\left( l(t) - \epsilon, u(t) + \epsilon \right).$
\end{lem}

\begin{proof}
Suppose $|a| \geq u(t) + \epsilon$.
%= \frac{2}{|t|^2} + \frac{1}{|t|} + \frac{2}{|t|^2}\sqrt{1+|t|} + \epsilon$. 
If we examine the modulus of the critical values we see that $|ta \pm 2\sqrt{a}| \geq |ta| - 2\sqrt{|a|} = |a|^{1/2}(|t||a|^{1/2} - 2) \geq (u(t)+\epsilon)^{1/2}(|t|(u(t)+\epsilon)^{1/2}-2)$.
Thus we have
\begin{align*}
& |ta \pm 2\sqrt{a}| \geq |ta| - 2\sqrt{|a|} \geq \\
&  \frac{2}{|t|} + 1 + \frac{2}{|t|}\sqrt{1+|t|} + \epsilon|t| - 2\sqrt{\frac{2}{|t|^2} + \frac{1}{|t|} + \frac{2}{|t|^2}\sqrt{1+|t|} + \epsilon}.
\end{align*}
We will rewrite two of the terms above: 
\begin{align*}
& \frac{2}{|t|} + \frac{2}{|t|}\sqrt{1+|t|} 
= \sqrt{\left(\frac{2}{|t|} + \frac{2}{|t|}\sqrt{1+|t|}\right)^2}
\\ & = \sqrt{\frac{4}{|t|^2} + \frac{4}{|t|^2}(1+|t|) + \frac{8}{|t|^2}\sqrt{1+|t|}}\\
& = \sqrt{\frac{8}{|t|^2} + \frac{4}{|t|} + \frac{8}{|t|^2}\sqrt{1+|t|}}
= 2\sqrt{\frac{2}{|t|^2} + \frac{1}{|t|} + \frac{2}{|t|^2}\sqrt{1+|t|}}.
\end{align*}
So,
\begin{align*}
& | ta \pm 2\sqrt{a}| \\
&\geq \frac{2}{|t|} + 1 + \frac{2}{|t|}\sqrt{1+|t|} + \epsilon|t| - 2\sqrt{\frac{2}{|t|^2} + \frac{1}{|t|} + \frac{2}{|t|^2}\sqrt{1+|t|} + \epsilon}\\
% &= 1 + \epsilon|t| + \frac{2}{|t|} + \frac{2}{|t|}\sqrt{1+|t|} - 2\sqrt{\frac{2}{|t|^2} + \frac{1}{|t|} + \frac{2}{|t|^2}\sqrt{1+|t|} + \epsilon}\\
&= 1 + \epsilon|t| + 2\sqrt{\frac{2}{|t|^2} + \frac{1}{|t|} + \frac{2}{|t|^2}\sqrt{1+|t|}} - 2\sqrt{\frac{2}{|t|^2} + \frac{1}{|t|} + \frac{2}{|t|^2}\sqrt{1+|t|} + \epsilon}.
\\ & = \colon 1+ g_t(\epsilon)
\end{align*}

We define the function $g_t(\epsilon)$ so that $| ta \pm 2\sqrt{a}| =    1+ g_t(\epsilon)$. Then to show  $|ta \pm 2\sqrt{a}| > 1 + \delta_1$ for some $\delta_1 > 0$, we'll show $g_t(\epsilon)>0$. 
First, note $g_t(0)= 0.$
Next,
$$
g_t'(\epsilon) = |t| - \frac{1}{\sqrt{\frac{2}{|t|^2} + \frac{1}{|t|} + \frac{2}{|t|^2}\sqrt{1 + |t|} + \epsilon}}.
$$
Now we observe:
\begin{align*}
&g_t'(\epsilon)>0 \Rightarrow |t| - \frac{1}{\sqrt{\frac{2}{|t|^2} + \frac{1}{|t|} + \frac{2}{|t|^2}\sqrt{1 + |t|} + \epsilon}} > 0 \\
& \Rightarrow  \ |t|\sqrt{\frac{2}{|t|^2} + \frac{1}{|t|} + \frac{2}{|t|^2}\sqrt{1 + |t|} + \epsilon} > 1
\end{align*}
$$ \Rightarrow |t|^2 \left( \frac{2}{|t|^2} + \frac{1}{|t|} + \frac{2}{|t|^2}\sqrt{1 + |t|} + \epsilon \right) > 1$$
$$\Rightarrow  2 + |t| + 2\sqrt{1 + |t|} + \epsilon|t|^2 > 1$$
$$\Rightarrow  1 + |t| + 2\sqrt{1 + |t|} + \epsilon|t|^2 > 0,$$
which is true for any $\epsilon > 0$ and any choice of $t$. 
Therefore, $g_t(\epsilon)$ 
is increasing with respect to $\epsilon$. Since $g_t(0) = 0$, if  $\epsilon > 0$, $g_t(\epsilon)>0$. 
Thus, there is a $\delta_1 > 0$ for which 
$1+g_t(\epsilon) > \delta_1.$
Therefore
$ | ta \pm 2\sqrt{a}| > 1 + \delta_1.$

By Corollary~\ref{cor:BS1} this means that for any choice of $a$ where $|a|\geq u(t) + \epsilon$ there exists an $N_1 \geq 2$ such that for all $n \geq N_1$ we have $ta \pm 2\sqrt{a} \not\in K(\Rnta)$. Thus whenever $|a| \geq u(t)+\epsilon$ there exists an $N_1 \geq 2$ such that for all $n \geq N_1$ we get $a \not\in M_n(\Rnta)$. We next use a similar process to show $a \not\in M_n(\Rnta)$ for $a$ to small to be in the annulus.

Suppose $|a| \leq l(t) - \epsilon$.
Then as above:
\begin{align*}
|ta \pm 2\sqrt{a}| &\leq |ta| + 2\sqrt{|a|}\\
&\leq \frac{2}{|t|} + 1 - \frac{2}{|t|}\sqrt{1+|t|} - \epsilon|t| + 2\sqrt{\frac{2}{|t|^2} + \frac{1}{|t|} - \frac{2}{|t|^2}\sqrt{1+|t|} - \epsilon}.
\end{align*}
In a similar calculation to the upper boundary case, looking at two terms from this expression, we can write
\begin{align*}
& \frac{2}{|t|} - \frac{2}{|t|}\sqrt{1+|t|} = -\sqrt{\left(\frac{2}{|t|} - \frac{2}{|t|}\sqrt{1+|t|}\right)^2}
 = -2\sqrt{\frac{2}{|t|^2} + \frac{1}{|t|} - \frac{2}{|t|^2}\sqrt{1+|t|}}.
\end{align*}
Note that $\frac{2}{|t|} - \frac{2}{|t|}\sqrt{1+|t|}$ is negative, which is why we choose the negative square root above. Using this we have
\begin{align*}
& |ta \pm 2\sqrt{a}| \\
&\leq \frac{2}{|t|} + 1 - \frac{2}{|t|}\sqrt{1+|t|} - \epsilon|t| + 2\sqrt{\frac{2}{|t|^2} + \frac{1}{|t|} - \frac{2}{|t|^2}\sqrt{1+|t|} - \epsilon}\\
&= 1 - \epsilon|t| + \frac{2}{|t|} - \frac{2}{|t|}\sqrt{1+|t|} + 2\sqrt{\frac{2}{|t|^2} + \frac{1}{|t|} - \frac{2}{|t|^2}\sqrt{1+|t|} - \epsilon}\\
&= 1 - \epsilon|t| - 2\sqrt{\frac{2}{|t|^2} + \frac{1}{|t|} - \frac{2}{|t|^2}\sqrt{1+|t|}} + 2\sqrt{\frac{2}{|t|^2} + \frac{1}{|t|} - \frac{2}{|t|^2}\sqrt{1+|t|} - \epsilon}.
\end{align*}
Note the final two root terms above differ only by an epsilon under one of the radicals, so their difference is a small negative number with absolute value less than epsilon. Hence this final expression can be written as $1$ minus $\epsilon|t|$ minus some small number. Thus there is a $\delta_2 > 0$ s.t.\  
$|ta \pm 2\sqrt{a}| < 1 - \delta_2.$\\
By Corollary~\ref{cor:BS1} this means that for any choice of $a$ where $|a|\leq l(t) - \epsilon$ there exists an $N_2 \geq 2$ such that for all $n \geq N_2$ we have $ta \pm 2\sqrt{a} \not\in K(\Rnta)$. Thus whenever $|a| \leq l(t) - \epsilon$ there exists an $N_2 \geq 2$ such that for all $n \geq N_2$ we get $a \not\in M_n(\Rnta)$.

Finally, let $N = \max(N_1, N_2)$ and we have our result; for all $n \geq N$ if $a \not\in \mathbb{A}(l(t)-\epsilon, u(t)+\epsilon)$ then $a \not\in M_n(\Rnta)$.
\end{proof}

We now have an annulus containing the boundedness locus. 
To refine, we show the boundedness locus lays in a neighborhood of $\mathcal{S}_t$.

\begin{proof}[Proof of Theorem~\ref{thm:Main_Theorem_C=ta}]
$\textbf{Notation:}$ We denote the $\epsilon$-neighborhood of a set $S$ by $N_\epsilon(S)$, and the ball of radius $\epsilon$ centered at a point $z$ by $B_\epsilon(z)$.

This proof generalizes the proof of Lemma 5.12 in \cite{BoydSchulz}. First, we apply Lemma~\ref{lem:MandelbrotInAnnulus} to get an $N_1$ such that for all $n\geq N_1$ we have $M_n(\Rnta) \subset \mathbb{A}(l(t)-\epsilon,u(t)+\epsilon)$. Note that the proof of the lemma implies that $\mathcal{S}_t \subset \mathbb{A}(l(t),u(t))$. We want to show that if $a \not\in N_\epsilon(\mathcal{S}_t)$ then $a \not\in M_n(\Rnta)$, so there are two cases to consider.

First, if $a \not\in \mathbb{A}(l(t)-\epsilon,u(t)+\epsilon)$, then by Lemma \ref{lem:MandelbrotInAnnulus} we have $a \not\in M_n(\Rnta)$.

Second, consider the set $H_{t,\epsilon} = \mathbb{A}(l(t)-\epsilon,u(t)+\epsilon) \setminus N_\epsilon(\mathcal{S}_t)$. This set is compact, so if we create an open cover of the set, $\{B_{\epsilon/2}(z) | z \in H_{t,\epsilon}\}$, it will have a finite subcover $\{B^1_{\epsilon/2},B^2_{\epsilon/2},...,B^m_{\epsilon/2}\}$.

Since $0 \not\in \mathbb{A}(l(t),u(t))$, we can take $\epsilon$ sufficiently small so that none of the $\{B^1_{\epsilon/2},B^2_{\epsilon/2},...,B^m_{\epsilon/2}\}$ include $0$. Therefore the square root function turns each $\{B^1_{\epsilon/2},B^2_{\epsilon/2},...,B^m_{\epsilon/2}\}$ into two open neighborhoods, and $\left\lbrace ta \pm 2\sqrt{B^k_{\epsilon/2}}\right\rbrace$ is a collection of open neighborhoods with each $\overline{B^k_{\epsilon/2}} \subset \CC\setminus N_{\epsilon/2}(\mathcal{S}_t)$. Now define $\delta_k$, for each $k$, and $\delta$ by
$$\delta_k = \min_{z\in \overline{B^k_{\epsilon/2}}} \left| \left| tz \pm 2\sqrt{z}\right| - 1 \right|, 
\text{ and }
\delta = \displaystyle\min_{1\leq k \leq m} \delta_k.$$

For this $\delta$ we apply Corollary~\ref{cor:BS1} to get an $N \geq N_1$ such that for all $n \geq N$ we have $K(\Rnta) \subset N_\delta(S^1)$. However, for each $a \in H_{t,\epsilon}$ we also have $a$ in one of the $B^k_{\epsilon/2}$, and therefore $||ta \pm 2\sqrt{a}|-1| > \delta$. This means that both $ta \pm 2\sqrt{a} \not\in K(\Rnta)$ and therefore $a \not\in M_n(\Rnta)$.

Thus we have that, for all $n \geq N$, if $a \not\in N_{\epsilon}(\mathcal{S}_t)$ then $a \not\in M_n(\Rnta)$. Therefore $M_n(\Rnta) \subset N_\epsilon(\mathcal{S}_t)$.
\end{proof}

\subsection{Future work}

For $\Rnta$ we identified a spine around which the boundedness locus must be centered, as well as a neighborhood around this spine in which the boundedness locus must lie. The next steps would be to find centers for where the baby $\mandel$'s should be, and then to prove explicitly that the baby $\mandel$'s are present.
 
To locate the centers, for each fixed $n\geq 2,$ Lemma~\ref{cor:BS2} was a result in the proof of the final proposition of \cite{BoydSchulz}. For the map $R_{n,a,c}$, for any $c \in \CC^*, n \geq 2$ this proof locates regions containing the distinct values $w_{n,k} = w_k$ such that for $a_k = w^{2n}_k,$ we have $v_+$ or $v_-$ equal to one of the roots $a^{1/2n}$ (for $n$ values of $k$ if $|c|\geq 1$ or $n-1$ values if $|c| < 1$). 

For $\Rnta$, $c$ is restricted to $c=ta$ for some $t\in \CC^*$. One would want to examine the proof of Proposition 5 in \cite{BoydSchulz} to add this restriction, to try to establish how many fixed critical points, thus potential centers of potential baby $\mandel$s, there are in this plane for each $t\neq 0.$ Looking at this proof, setting $t=ca$ means rather than solving 
$
w^n = \frac{w}{2} - \frac{c}{2},
$
for $w$ with $w^{2n}=a$, one must look for (distinct) solutions of
$
w^n = \frac{w}{2} - \frac{t w^{2n}}{2},
$
or
$
\frac{w^n-1}{2w^{2n-1}} -t =0.  
$

%%%%%%%%%%%%%%%%%%%%%%%

\end{document}